\documentstyle[12pt]{article}
\renewcommand{\baselinestretch}{1.3}

\newenvironment {qth}[1]{{ }\\[2ex] {\bf Theorem #1:} \it}{{ }\\[2ex]}

\newtheorem {th}{Theorem}[section]
\newtheorem {lem}[th]{Lemma}
\newtheorem {pr}[th]{Proposition}
\newtheorem {cor}[th]{Corollary}

\newtheorem{conj}{Conjecture}
\def\Cox{\hfill \Box}

\def\dd{\delta}

\def\iid{\mbox{independent, identically distributed}}

\def\ee{\epsilon}
\def\E{{\bf{E}}}
\def\P{{\bf{P}}}
\def\N{\hbox{I\kern-.2em\hbox{N}}}
\def\R{\hbox{I\kern-.2em\hbox{R}}}
\def\Z{{\bf{Z}}}
\def\B{{\cal{B}}}

\def\|{\, | \, }
\def\one{{\bf 1}}
\def\0{\hat{0}}
\def\1{\hat{1}}
\def\rr{{\bf r}}
\def\vw{{\overline {vw}}}
\newcommand{\lc}{\left\lceil}
\newcommand{\lf}{\left\lfloor}
\newcommand{\rc}{\right\rceil}
\newcommand{\rf}{\right\rfloor}

\begin{document}

\begin{center}
{\large \bf PERCOLATION, FIRST-PASSAGE
PERCOLATION, AND COVERING TIMES FOR
RICHARDSON'S MODEL ON THE $n$-CUBE} \\
\vspace{1ex}
[Short title: PERCOLATION ON THE CUBE] \\
\vspace{2ex}
{\sc By James Allen Fill and Robin Pemantle}\\
\vspace{1ex}
%%%%I've updated my affiliation on the two cover pages and at the end.
{\em The Johns Hopkins University and University of Wisconsin-Madison\/}
\end{center}
\newpage

\begin{center}
{\large \bf PERCOLATION, FIRST-PASSAGE
PERCOLATION, AND COVERING TIMES FOR
RICHARDSON'S MODEL ON THE $n$-CUBE} \\
\vspace{1ex}
[Short title: PERCOLATION ON THE CUBE] \\
\vspace{2ex}
{\sc By James Allen Fill} \footnote{Research supported by the National Security
Agency under Grant Number MDA904-89-H-2051.}
{\sc and Robin Pemantle} \footnote{Research supported by a National
Science Foundation postdoctoral fellowship and by a Mathematical
Sciences Institute postdoctoral fellowship.}\\
\vspace{1ex}
{\em The Johns Hopkins University and University of Wisconsin-Madison\/}

\footnotetext[3]{{\em AMS\/} 1991 {\em subject classifications\/}.  Primary
60K35; secondary 60C05.}
%**%60K35 refers to percolation theory.  60C05 is Combinatorial probability.


\end{center}

\vspace{2ex}

{\bf ABSTRACT:} \break
Percolation with edge-passage probability $p$ and
first-passage percolation
are studied for the $n$-cube $\B_n = \{ 0 , 1 \}^ n$ with
nearest neighbor edges.  For oriented and unoriented percolation, $p = e/n$
and $p = 1/n$ are the respective critical probabilities.  For oriented
first-passage percolation with i.i.d.\ edge-passage times having a density of
$1$ near the origin, the percolation time (time to reach the opposite corner
of the cube) converges in probability to $1$ as $n \rightarrow \infty$.  This
resolves a conjecture of David Aldous.
When the edge-passage distribution is standard exponential,
the (smaller) percolation time
for unoriented edges is at least $0.88$.

These results are applied to Richardson's model on the (unoriented) $n$-cube.
Richardson's model, otherwise known as the contact process
with no recoveries, models the spread of infection as a
Poisson process on each edge connecting an infected node to
an uninfected one.  It is shown that the time to cover the
entire $n$-cube is bounded between $1.41$ and $14.05$ in probability as $n
\rightarrow \infty$.

\section{Introduction and notation}

Percolation theory, broadly speaking, is the study of connectivity
in a random medium.  Grimmett (1989) gives the following example
of a question addressed by percolation theory.  ``Suppose we
immerse a large porous stone in water.  What is the probability
that the centre of such a stone is wetted?''  The mathematical
model for this is graph-theoretic.  Small volume elements of
the stone become vertices of the integer lattice (of dimension
three in this case).  Neighboring volume elements may be connected
by a channel broad enough to allow the flow of water, or
they may not be.  Model this by letting each pair of
neighboring vertices independently have a connecting edge
with probability $p$, for some parameter $p$.  The center of
the stone is then wetted if and only if the connected component
of this random subgraph containing the center extends to
the surface.

Oriented percolation is a variant on this, where each edge has a
particular orientation, and the water may pass only in that
direction if at all.  For example, the stone may be subjected
to water only from above, with flow of water from one volume
element to a neighboring one occurring (due to gravity) only
when the latter is lower.  First-passage percolation is a
similar model, the difference being that the passage of
water through a channel is not simply a yes or no event, but
takes an amount of time depending on the breadth of the channel.
Thus each edge between neighboring vertices, instead
of being randomly included or excluded, is assigned a
random (i.i.d.) passage time.  The question is not whether,
but when the center of the stone first gets wet.

Richardson's model is a stochastic process on a graph that
begins with one vertex {\em infected\/} and evolves by
transmission of the infection according to an i.i.d.\ Poisson
process on each edge: if the Poisson process on the edge
connecting $v$ to $w$ has a point of increase at time $t$
and one of $v$ or $w$ was infected before time $t$, then
both are infected after time $t$.  Questions about Richardson's
model can be reduced to questions about first-passage
percolation with i.i.d.\ exponentially distributed passage
times.  In Pemantle (1988), questions about recurrence
or transience of random walks in random environments on
trees are essentially reduced to oriented first-passage
percolation on trees.  Lyons (1990) connects random walks
and electrical networks with percolation on trees.  These
equivalences (along with the desire to generalize whenever
possible) are among the reasons to study percolation on graphs such
as trees and $n$-cubes which do not model any actual stones.

Classical percolation theory, i.e., on $\Z^ d$, is geometric.
The main arguments---counting countours, piecing together
sponge crossings, viewing the process from its left edge---are
all pictorial.
Infinite binary trees were introduced as a way to get graphs
which were in some
sense limits of $\Z^ d$ as $d \rightarrow \infty$.  Trees are
in general easier to study than integer lattices.  For
example, the first-passage percolation problem for binary trees
is essentially a large deviation calculation: if every path
of length $n$ in the tree were disjoint then the passage
time would be the minimum of $2^ n$ different sums of $n$
i.i.d.\ random variables; such a calculation is quantitative,
requiring no picture.  The argument in Pemantle (1988)
consists mostly of showing that the slight overlapping
of paths is inconsequential.  Unfortunately, trees
are locally quite different from lattices (despite being
called {\em Bethe lattices\/}) and are therefore not very
satisfactory in the role of limiting lattices.
The $n$-cube is an alternative
way to capture the high-dimensional limiting behavior of
integer lattices without altering the local connectivity
properties.  The solutions presented here to percolation
problems on $n$-cubes use large deviation and second
moment estimates and are thus closer to the tree case
than to the integer lattice case.  In this respect, the method
is very similar to the arguments used by Cox and Durrett (1983)
for the analogous problem on $(\Z^ + )^ d$ in high
dimensions.  The results
resolve affirmatively conjecture G7.1 of Aldous (1989), despite
a remark there indicating that second moment methods shouldn't
work.  (To be fair: the second moment method bounds the probability
in question away from $0$ but a variance reduction trick is needed
to get it equal to $1$.)

Notation will be as follows.  Let $\B_n$ be the Boolean algebra
of rank $n$ whose elements are ordered $n$-tuples of
$0$'s and $1$'s.  It has a bottom element $\0 = (0, \ldots , 0)$
and a top element $\1 = (1, \ldots , 1)$.  It is also useful to
regard elements of $\B_n$ as subsets of $\{ 1 , \ldots n\}$,
where the set $A$ corresponds to the sequence with a $1$
in position $i$ if and only if $i \in A$.
The collection of all $A$ of a given cardinality $k$ is called the
$k^ {th}$ {\em level\/} of $\B_n$.
The $n$-dimensional cube, or $n$-cube, is the graph whose vertices
are the elements of $\B_n$ and whose edges connect each set $A$
to $A \cup \{j\}$ for each $j \notin A$.  Sometimes we
must think of the edges as
oriented from $A$ to $A \cup \{j\}$.  Representing the vertices
of the $n$-cube as the standard basis in $\R^ n$ makes each such
edge parallel to the unit vector $e(j)$ connecting $\0$
to $(0, \ldots , 0 , 1 , 0 , \ldots , 0)$ with the $1$ in
position $j$.  Then a path from $\0$ to $\1$ of length $n$
can be represented as a permutation $(\pi (1) , \ldots ,
\pi (n))$ where the $k$th edge in the path is parallel
to $e (\pi (k))$ and in fact connects the set $A = \{\pi (i) :
i < k \}$ to the set $A \cup \{\pi (k)\}$.  Similarly, a path
connecting a set $A$ to a set $B \supseteq A$ may be represented as
a permutation of $B \setminus A$.

Let $X_{vw}$ be i.i.d.\ random variables as $\vw$ ranges
over all oriented edges of the $n$-cube.
The {\em first-passage time\/} from $\0$ to $\1$, or {\em percolation time\/},
is defined as the minimum over paths $\gamma = (\0 , v_1 , v_2 , \ldots ,
v_{n-1}, \1)$ from $\0$ to $\1$ of $T_n (\gamma )$, where $T_n
(\gamma)$ is
the sum of $X_{v_i , v_{i+1}}$ along the edges of $\gamma$.
In the application of first-passage
percolation to Richardson's model, the common distribution
of the $X_{vw}$'s is exponential with mean $1$.  For our basic
first-passage
percolation result itself, the assumption of exponential edge-passage times
simplifies the derivation of certain large deviation estimates but is not
necessary; we state the theorem
for more general distributions of the passage time,
though probably still greater
generality is possible. For the integer lattice in high dimensions, Kesten
(1984) gives sufficient conditions on the behavior of the
common density of the edge-passage times near the origin for similar
calculations to work.
It is unlikely that the same conditions are sufficient for our problem;
at any rate, our paper does not address this issue.

Of course, for ordinary percolation, the common distribution
of the $X$'s takes on the two values ``open''
and ``closed'' with respective probabilities $p$ and $1-p$,
for some $p$.  The basic question then is to compute the probability that there
is a path from $\0$ to $\1$ consisting only of open edges.

The rest of the paper is organized as follows.  Section~2 contains
a description of the second moment method, and its enhancement by a variance
reduction technique, that is used to analyze oriented percolation and oriented
first-passage percolation.
Section~2 also presents lemmas that count oriented paths in
the $n$-cube.

%%%% Here begins some revision integrating statements of the main theorems
%%%% into the introductory discussion.
Section~3 discusses oriented percolation and oriented first-passage
percolation.  The probability of oriented percolation approaches a limit
when $n \rightarrow \infty$ and $np$ is constant, as given in the following
theorem.
\begin{qth}{3.2}
Let each edge of $\B_n$ be independently open with probability
$p=c/n$.  Then
$\P(\0 \mbox{\rm \ is connected to $\1$ by an oriented open path})$
converges to a limit as $n \rightarrow \infty$.
The limit is $0$ if $c < e$ and is $(1-x(c))^ 2$ if $c \geq e$, where
$x(c)$ is
the extinction probability for a Poisson($c$) Galton--Watson process, namely,
the solution in $(0,1)$ to $x=e^ {c(x-1)}$.
\end{qth}
For the oriented first-passage time we have
\begin{qth}{3.5} 
Let the edges of $\B_n$ be assigned i.i.d.\ positive random
passage times with common density $f$, and assume that
$|f(x) - 1| \leq K x$ for some $K$ and all $x \geq 0$.  Then
the oriented first-passage percolation time $T = T^ {(n)}$ for $\B_n$
converges
to $1$ in probability as $n \rightarrow \infty$.
\end{qth}

Section~4 treats unoriented percolation, for which the critical probability is
shown to be $1/n$.  The result here, which may already be part of the
percolation folklore, is
\begin{qth}{4.1}
Let each edge of $\B_n$ be independently open with probability $p = c/n$,
$0 <
c < \infty$.  Then $\P(\0$ {\rm is connected to $\1$ by an (unoriented) open
path)}$\to (1 - x(c))^ 2$, where $x(c)$ is, as in Theorem~\ref{ordinary},
the
extinction probability for a Poisson($c$) Galton--Watson process.
\end{qth}

Section~5 gives an argument of Durrett (personal communication)
provideing a lower bound for the unoriented first-passage time
when the common distribution of the passage times is exponential
by comparing the process to a branching translation process (BTP), which
is similar to a branching random walk.  
\begin{qth}{5.3 (Durrett)}
As $n \rightarrow \infty$, the time $\tau_n$ of first population
of $\1$ in BTP converges in probability to $\ln(1+\sqrt{2})
\doteq 0.88$.  Consequently, $\P (T_n \leq \ln(1+\sqrt{2}) - \ee)
\rightarrow 0$.
\end{qth}
Since the oriented
first-passage time is an upper bound for the unoriented
first-passage time, the unoriented first-passage time is
thereby bounded as $n \rightarrow \infty$
between two fairly close constants, namely, $0.88$ and $1$.  We remark
that for oriented first-passage times, lower bounds are easy (first 
moment calcluation) while upper bounds equalling these lower bounds 
are more difficult (second moment estimates).
For unoriented first-passage times, we do not know
how to bridge the gap between the lower and upper bounds.

Finally, the sixth section discusses the cover time in Richardson's
model, i.e. the first time that all sites are infected.  As $n \rightarrow
\infty$, the cover time is bounded in probability by a constant; this
is shown by improving Theorem~3.5 so as to get an exponentially small 
bound on the probability of a vertex failing to be reached by a time 
$c = 3+2\ln (4+2\sqrt{3}) \doteq 7.02$.  This is certainly not sharp
though it improves on the previous best upper bound for the covering
time which was of order $\ln n$.  A lower bound in probability of
$\ln (2) + {1 \over 2} \ln (2 + \sqrt{5})$ is also given, showing that
the covering time is separated from the single-vertex first-passage
time of 1.  
\begin{qth}{6.4 and Corollary 6.3}
For any $\ee > 0$, $\P (A (2c+\ee ) = \B_n ) \rightarrow 1$
as $n \rightarrow \infty$, where $c = 3+2\ln (4+2\sqrt{3}) \doteq 7.02$.
On the other hand, for any
$\ee > 0$, $\P (A(\frac{1}{2} \ln(2 + \sqrt{5}) + \ln 2 - \ee) = \B_n)
\rightarrow 0$ as $n \rightarrow \infty$.
\end{qth}

A future paper will address the discrete-time analogues of these problems,
which are closely related to the so-called broadcasting problem discussed by
Feige, et al.~(1990) and others.

\section{Preliminaries: the second moment method and path counting}

\subsection{The second moment method, with variance reduction}

Aldous (1989, Lemma A15.1) gives the following lemma as the cornerstone
of the so-called second moment method.  The proof is a simple application of
the Cauchy--Schwarz inequality.

\begin{lem}[Second moment method] \label{2nd moment method}
Let $N$ be a nonnegative real random variable with $\E N^ 2
< \infty$.  Then $\P (N > 0) \geq (\E N)^ 2 / \E N^ 2$.    $\Cox$
\end{lem}

In our applications, the random variable $N \equiv N_n$ will be indexed by
the
dimension $n$ of the cube under consideration. When the variance of $N$ is
$o((\E N)^ 2)$ as $n$ tends to infinity, the inequality shows that
$\P(N > 0)
= 1 - o(1)$.  On the other hand, if we can only show the variance of $N$ to be
$O((\E N)^ 2)$, then the conclusion is weaker, namely, that $\P (N > 0)$ is
bounded away from $0$.  One of the purposes of the present work is to point out
that the weaker conclusion may often be strengthened to the former by a simple
variance reduction trick.  Consider an auxiliary random variable $Z$ which
absorbs most of the variance of $N$, in the sense that the
conditional second moment $\E (N^ 2 \| Z)$ is
only $(1+o(1))$ times the conditional squared first moment
$(\E (N \| Z))^ 2$ uniformly over a set of values of $Z$ of probability
$1 -
o(1)$. Then uniformly for those values of $Z$, $\P (N > 0 \| Z)$ is
$1 - o(1)$,
and hence $\P(N > 0)$ is also $1 - o(1)$.

To apply the second moment method to oriented percolation on the
$n$-cube, let $N$ be the number of paths of length $n$
from $\0$ to $\1$ consisting entirely of open edges.  Then
$\0$ is connected to $\1$ by an oriented path of open edges
if and only if $N > 0$.
For large enough values of $p = \P (\mbox{edge is open})$, the (unenhanced)
second moment method will get $\P(N > 0)$ bounded away from $0$ by
showing that $(\E N^ 2) / (\E N)^ 2$ is bounded, in the following manner.
A similar argument, to be detailed later, using the enhanced second moment
method will yield the exact limiting value of the percolation probability $\P
(N > 0)$.

Since $\E N$ is the sum over paths $\gamma$ of the probability that
$\gamma$ consists entirely of open edges and $\E N^ 2$ is the sum over
pairs of paths $(\gamma , \gamma')$ of the probability
that both paths consist entirely of open edges, the quotient $\E N^ 2 / (\E
N)^ 2$ would be precisely $1$ if the events that $\gamma$ consists
entirely of open edges and that $\gamma'$ consists
entirely of open edges were always independent.  Of course
they are not independent when $\gamma = \gamma'$, but also
they are not independent when $\gamma$ and $\gamma'$ have
any edges in common, and the covariance of their indicators is greater
the more edges that $\gamma$ and $\gamma'$ share.  Thus the argument rests
on showing that the number of pairs $(\gamma , \gamma' )$ sharing
a lot of edges is small.

The oriented first-passage percolation
problem is handled similarly, with $N$ being the number of paths whose
total passage time is at most $M$; if $N > 0$ with high probability,
then the passage time is less than $M$ with high probability.

\subsection{Counting oriented paths}

Section~2.2 contains the path-counting lemmas
needed to execute the second moment method.
By symmetry, it will only be necessary to consider the case
where $\gamma$ is the path given by $(1, \ldots , n)$ in
the permutation representation.  Let $f(n,k)$ be the
number of paths $\gamma'$ that share precisely $k$ edges
with $\gamma$.  Let $F(n,k) = \sum_{n \geq l \geq k} f(n,l)$.
Most of the time it will suffice to bound $F(n,k)$ and
observe that $f(n,k) \leq F(n,k)$.

\begin{lem} \label{small k}
Let $K(n)$ be any function that is $o(n)$ as
$n \rightarrow \infty$.  Then $f(n,k) \leq F(n,k)
\leq (1 + o(1)) (k+1) (n-k)!$ as $n \rightarrow \infty$ uniformly in
$k$ for $k \leq K(n)$.  Furthermore, for $k>0$,
consider paths that agree with $\gamma$ only in an initial segment
and a final segment;
it is these paths that matter for $F(n,k)$, in the sense that
all the rest of the paths only
contribute $o((k+1)(n-k)!)$, uniformly in $0 < k \leq K(n)$.
\end{lem}

\begin{lem} \label{large k}
Suppose $k \geq n - n^ {3/4}/2$.  Then $f(n,k) \leq F(n,k) \leq
(n-k+1)(2n^ {7/8})^ {n-k}$.
\end{lem}

\begin{lem} \label{middle k}
Suppose $k \leq n - 5e (n+3)^ {2/3}$.  Then, for $n
\geq 25$, $f(n,k) \leq n^ 6 (n-k)!$, and, for
$n$ so large that $\lc 5 e (n+3)^ {2/3} \rc \leq \lc n^ {3/4} / 2 \rc$,
$$
F(n,k) \leq 2n^ 6 (n-k)! + \lc 5 e (n+3)^ {2/3} \rc
(2n^ {7/8})^ {\lc 5 e (n+3)^ {2/3} \rc - 1}.
$$
\end{lem}

The following lemma will be needed when we apply the enhanced second moment
method.

\begin{lem} \label{Fsub1}
Consider now a Boolean lattice $\B_{n+2L}$ of size $n + 2L$
for some positive integer $L$.  Let $x_1$ and $x_2$ be distinct
vertices in level $L$ of $\B_{n+2L}$
and let $y_1$ and $y_2$ be distinct vertices in level $n+L$
of $\B_{n+2L}$.  Assume that $x_i$ lies below $y_i$ for $i=1,2$.
Let $H (n,L,k,x_1,x_2,y_1,y_2)$ be the maximum over paths $\gamma$
from $x_1$ to $y_1$ of the number of paths from $x_2$ to
$y_2$ that have at least $k$ edges in common with $\gamma$.
Let $F_1 (n,L,k) = \max_{x_1,x_2,y_1,y_2}
H(n,L,k,x_1,x_2,y_1,y_2)$.
Then for fixed $L$ and any function $K(n)$ that is $o(n)$ as
$n \rightarrow \infty$, the following hold:
\begin{equation} \label{small Fsub1}
F_1 (n,L,k) = o((k+1)(n-k)!) \mbox{ uniformly in } 0 < k \leq K(n);
\end{equation}
\begin{equation} \label{large Fsub1}
F_1 (n,L,k) \leq (n-k+1) (2n^ {7/8})^ {(n-k)} \mbox{ for } k
\geq n - n^ {3/4}/2;
\end{equation}
\begin{equation} \label{middle Fsub1}
F_1 (n,L,k) \leq 2n^ 6 (n-k)! +  \lc 5 e (n+3)^ {2/3} \rc
(2n^ {7/8})^ {\lc 5 e (n+3)^ {2/3} \rc - 1} \mbox{ for } k \leq n -
5e(n+3)^ {2/3}.
\end{equation}
For the last inequality we require $\lc 5 e (n+3)^ {2/3} \rc \leq \lc
n^ {3/4} / 2 \rc$.
\end{lem}

The following notation is common to the proofs of Lemmas~\ref{small
k}~--~\ref{middle k}.
If $\gamma'$ has precisely $k$ edges in common with $\gamma$
and is given by the permutation $\pi$, then let $r_1 , \ldots ,
r_k$ be the
terminal positions of the shared edges.  In other words, let $r_1$ be
minimal
so that the set $\{ 1 , \ldots , r_1 -1 \}$ is equal to the
set $\{ \pi (1) , \ldots , \pi (r_1 - 1) \}$ and $\pi (r_1)
= r_1$; let $r_2$ be the next such value, and so on.  By convention,
let $r_0$
always be $0$ and $r_{k+1}$ always
be $n+1$.  Write
$\rr = \rr (\gamma' )$ for the sequence $(r_0, \ldots , r_{k+1})$.
Write $s_i = r_{i+1} - r_i$, $i = 0, \ldots, k$.  For any sequence
$\rr_0 = (r_0 , \ldots , r_{k+1})$ with $0 = r_0 < r_1 <
\cdots < r_k
< r_{k+1} = n+1$, let $C(\rr_0 )$ denote the number of paths
$\gamma'$ with $\rr (\gamma' ) = \rr_0$.  Then it is easy to
see that $C(\rr ) \leq G(\rr )$, where
\begin{equation} \label{Gdef}
G(\rr ) = \prod_{i = 0}^ k (s_i - 1)! ,
\end{equation}
since the values $\pi (r_i + 1) , \ldots , \pi (r_i + s_i - 1)$
must be a permutation of $\{ r_i + 1 , \ldots , r_i + s_i - 1 \}$.
For Lemma~\ref{middle k}, the more precise bound $C(\rr )
\leq G_1 (\rr )$ will be necessary, where $G_1$ is defined by
\begin{equation} \label{G1def}
G_1 (\rr ) = \prod_{i = 0}^ k [(s_i - 1)! - 1 + \dd\sb
{1,s_i}]
\end{equation}
and $\dd_{1 , s_i}$ is $1$ if $s_i = 1$ and $0$ otherwise.
To see this inequality, recall that $\gamma'$ must not have any
common edges with $\gamma$ strictly between $r_i$ and $r_{i+1}$, and
therefore,
for $s_i \neq 1$, at least one permutation of $\{ r_i + 1 , \ldots ,
r_{i+1} - 1
\}$ (namely, the identity permutation) is ruled out for the values of
$\pi (r_i + 1) , \ldots , \pi(r_{i+1} - 1)$.
Note that $s_i (\rr (\gamma'))$ can never be $2$, a fact which
is used liberally in the proofs of Lemmas~\ref{large k} and~\ref{middle k}.
Of course, $f(n,k)$ is equal to the sum of $C(\rr )$ over
all sequences $\rr$ satisfying $0 = r_0 < r_1 < \cdots < r_k
< r_{k+1} = n+1$.  Furthermore, $F(n,k)$ is at most the sum of $G(\rr )$
over all such sequences, since $G(\rr )$ counts the paths $\gamma'$
with $\rr \subseteq \rr (\gamma' )$.
It remains then to get bounds on these sums.  We use the following
facts about factorials.

\begin{pr} \label{factorials}
\begin{itemize}

\item[{\rm (i)}] Factorials are log-convex, i.e., $a! b! \leq (a+j)! (b-j)!$
for $a \geq b \geq j \geq 0$.

\item[{\rm (ii)}] $(a! - 1) ((a+j)! - 1) \leq ((a-1)! - 1) ((a+j+1)! - 1)$
as long as $a \geq 4$ or $a = 3$ and $j \geq 1$.

\item[{\rm (iii)}] $a > b > 0$ implies $(a! - 1) / (b! - 1) > a! / b!
> (a/e)^ {a-b}$.
\end{itemize}
\end{pr}

\noindent{\bf Proof:}  (i) follows from (in fact, is equivalent to) the fact
that
$(a+1)! / a! = a+1$ is increasing.  (ii) is easy to verify.
The first inequality of (iii) is trivial.  To prove the second one,
use Stirling's formula $a! e^ {-1/(12 a)} < a^ a e^ {-a}
\sqrt{2 \pi a} < a!$
to get
$$a!/(b! (a/e)^ {a-b}) > {a^ a e^ {-a} \sqrt{2 \pi a} \over
    b^ b e^ {-b} \sqrt{2 \pi b} e^ {1/(12b)}} (e/a)^ {a-b}
    = (a/b)^ {b+1/2} e^ {-1/(12 b)},$$
and taking the logarithm of the last expression gives at least
$$(b + 1/2) \ln(1 + 1/b) - 1/(12 b) > (b+ 1/2) (1/b - 1/b^ 2) -
1/(12b) > 0$$
for $b \geq 2$.  For $b = 1$, the result follows directly
from Stirling's formula.
$\Cox$

\noindent{\bf Proof} of Lemma~\ref{small k}:  For fixed $\rr$, let
$j = j(\rr ) = \max_i(s_i - 1)$.
Consider separately the cases $j \leq n-4k$ and $j > n-4k$.
The idea is that $G(\rr )$ is small in the first case, and while
$G$ is not so small in the second case, there are not too many
sequences with such large values of $j$.
In the first case, $G(\rr )$ must be no more than $(n-4k)!
(3k)!$.  To see this, note that by log-convexity of factorials,
the product in~(\ref{Gdef}) is maximized subject to $j \leq n-4k$
at some $\rr$ for which $j(\rr )$ is equal to $n-4k$, in which case
since $\sum_i (s_i -1) = n-k$, $G (\rr )$ is at most $[\max\sb
i(s_i - 1)]!
[\sum_i(s_i - 1) - \max_i(s_i -1)]! \leq (n-4k)!(3k)!$.
Meanwhile, the number of sequences $\rr$
under consideration is at most ${n \choose k}$.  Thus the total
contribution from this case is at most
\begin{eqnarray*}
&& (n-4k)! (3k)! n! / (k! (n-k)!) \\[2ex]
&  = & (n-k)! {(n-4k)! \over (n-k)!} {n! \over (n-k)!}
    {(3k)! \over k!} \\[2ex]
& = & (n-k)! {[n(n-1) \cdots (n-k+1)] [3k (3k-1) \cdots (k+1)]
    \over (n-k) (n-k-1) \cdots (n-4k+1)} \\[2ex]
& \leq & (n-k)! \left [ {(3k)^ 2 n \over (n-4k+1)^ 3} \right ]^ k
\end{eqnarray*}
The term inside the square brackets converges to $0$ as
$n \rightarrow \infty$ uniformly in $k$ as long as $k \leq K(n)
= o(n)$, so for large $n$ and $k \leq K(n)$, the contribution
from the case $j(\rr ) \leq n-4k$ is $o((n-k)!)$, uniformly in these $k$.

The second contribution is from the terms with $j > n-4k$.
Since factorials are log-convex, it follows (as in the preceding case) that
$G(\rr )$ can be no greater than $j(\rr )! (n-k-j(\rr ))!$.  The
number of sequences $\rr$ for which $j(\rr )$ is any fixed value
$j_0$ is at most $k+1$ times the number of sequences $\rr$
with $s_0 - 1 = j_0$, which is in turn at most $(n - j_0 - 1)!
/ ((k-1)! (n- j_0 - k)!)$.
So the contribution for fixed $j_0$ is at most
$${(k+1) j_0! (n-j_0-k)! (n-j_0 -1)! \over (k-1)! (n-j_0-k)!}
    = {(k+1) \over (k-1)!} j_0 ! (n - j_0 - 1)!.$$
Changing $j_0$ to $j_0+1$ multiplies this by $(j_0+1) / (n-j_0
-1)
> (n-4K(n)) / (4K(n))$, which tends to infinity by assumption on $K(n)$;
hence the sum is bounded by $(1 + o(1))$ times the single term with the maximum
value of $j_0$, namely, $j_0 = n-k$.  The contribution from
this term is at most $[(k+1) / (k-1)!] (n-k)! (k-1)! = (k+1)(n-k)!$.
The case $j(\rr ) = n-k$ covers precisely those sequences $\rr$
consisting of an initial segment and a final segment; thus both
statements in the lemma have been proved.   $\Cox$

\noindent{\bf Proof} of Lemma~\ref{large k}:  Let $m = n-k$ be
the number of positions in which the edges of $\gamma$ and
$\gamma'$ differ and let $j(\rr )$ this time be the
number of runs of consecutive positions that constitute these
$m$ edges.  For example, if $n=10$ and $k=6$ with $(r_0 ,
\ldots , r_{k+1}) = (0, 3, 4, 5, 8, 9, 10, 11)$
then $j(\rr ) = 2$ because the edges not shared are in
two runs of consecutive positions, namely, $1,2$ and $6,7$.
Now the inequality $F(n,k) \leq \sum_\rr G(\rr )$ can be
broken up as
$$
F(n,k) \leq \sum_l g(n,k,l) h(n,k,l),
$$
where $g(n,k,l)$ is an upper bound for $G(\rr )$ over
sequences $0 = r_0 < \cdots < r_{k+1} = n+1$ with
$j(\rr ) = l$ and $h(n,k,l)$ is the number of sequences
with $j(\rr ) = l$.

As remarked before, each such run of consecutive positions contains
at least two positions, since otherwise the single edge in the run
would also be a shared edge.  Thus $l \leq m/2$, and for fixed
$l$, the number of $\rr$ for which $j(\rr ) = l$ is equal to
the number of ways of choosing $m$ positions out of $n$ in $l$
clusters of size at least $2$ each.  Treating each cluster
as a unit, there are ${{n-m+l} \choose l}$ ways of locating the $l$
clusters among the $n-m$ shared edges.  Since each cluster
has at least two positions, there are $m-2l$ extra positions
to be distributed among the $l$ clusters.  The number of ways to
do this is the number of ways to drop $m-2l$ indistinguishable
balls into $l$ distinguishable boxes, which is ${{m-l-1} \choose {l-1}}$.
Thus $h(n,k,l) = {{n-m+l} \choose l} {{m-l-1} \choose {l-1}}$.  The first
factor is increasing in $l$, thus maximized when $l=\lf m/2 \rf$,
while the second factor is at most ${m \choose l}$, which is also
maximized at $l=\lf m/2 \rf$.  Thus, assuming for simplicity that $m$ is even,
$$
h (n,k,l) \leq {{n - m/2} \choose {m/2}} {m \choose {m/2}}.
$$

Using again the log-convexity of factorials, $G(\rr )$ is less than
or equal to the product of the factorials of the cluster sizes, which
is maximized when all clusters have size $2$ except for
one large cluster.  This gives $g(n,k,l) \leq 2^ {l-1} (m-2l+2)! \leq
2^ {m/2}
m!$, and thus
\begin{eqnarray*}
\sum_{l=0}^ {m/2} g(n,k,l) h(n,k,l) & \leq & ({m \over 2} + 1)
{{n - m/2} \choose
{m/2}}
    {m \choose {m/2}} 2^ {m/2} m! \\[2ex]
& \leq & ({m \over 2} + 1) {n^ {m/2} \over (m/2)!} {m^ {m/2} \over
(m/2)!}
    2^ {m/2} m! \\[2ex]
& = & ({m \over 2} + 1) (2mn)^ {m/2} {m \choose {m/2}} \\[2ex]
& \leq & ({m \over 2} + 1) (8mn)^ {m/2} .
\end{eqnarray*}
The same inequality $\sum_l g(n,k,l) h(n,k,l) \leq
({m \over 2} + 1) (8mn)^ {m/2}$ can be established similarly when $m$ is
odd.
Now the hypothesis that $m < n^ {3/4} / 2$ implies that the bound
$({m \over 2} + 1) (8mn)^ {m/2}$ is at
most $(m+1) (2n^ {7/8})^ m$.    $\Cox$

\noindent{\bf Proof} of Lemma~\ref{middle k}:  The idea this
time is to define a weight $w(\rr )$ with the property that the sum
of the $w(\rr )$ over all sequences $0 = r_0 < \cdots < r_{k+1} = n+1$
is less than $1$.  Then
$$
\sum_\rr C(\rr) = \sum_\rr w(\rr) [C(\rr) / w(\rr)]
\leq \max_\rr[C (\rr) / w (\rr)],
$$
which will be shown to be as small as required.
Let $n_i$ denote $n - r_i$ and let $s_i = r_{i+1} - r_i$ as
before,
viewing the sequence of $s_i$'s as a function of $\rr$.
The weight $w$ is defined by
$$
w(\rr) = \prod_{i=0}^ k Q(n_i, s_i),
$$
where
$$
Q(m , s) = \left \{ \begin{array} {lcl} (m + 1) / (m + 3) & : &
    s = 1 \\ 1/(m+3) & : & s = 3 \\ 1 / [(m+3)(m+2)] & : &
    s \geq 4 . \end{array} \right.
$$
Note that $Q(m,1) + \sum_{s=3}^ {m+1} Q(m,s) < 1$ for each fixed $m$;
hence by a
simple induction argument the sum of $w (\rr )$ over all sequences of a given
length $k$ is less than $1$.  It remains to bound for fixed $n$ and $k$ the
quantity $\max_\rr[C(\rr ) / w (\rr )]$.

The procedure will be to find for each $\rr$ a quantity $D(\rr)$ and a pair
$\rr^ {(1)}, \rr^ {(2)}$ such that
$C(\rr)/w(\rr) \leq \frac{1}{2} (n+3)^ 2 D(\rr)$, $D(\rr) \leq D(\rr\sp
{(1)}) \leq
(n+3)^ 2 D(\rr^ {(2)})$, and $D(\rr^ {(2)}) \leq D(\rr_0)$.  This
will give
$C(\rr)/w(\rr) \leq \frac{1}{2} (n+3)^ 4 D(\rr_0)$, and the calculation
$D(\rr_0) = [(n-k)!-1] (n-k+3) (n-k+2)$ will follow directly from the
definition of $D$.  Then
$$
\frac{C(\rr)}{w(\rr)} \leq \frac{1}{2} (n+3)^ 6 (n-k)! \leq n^ 6 (n-k)!
$$
for $n \geq 25$, as desired.

To begin, recall the bound~(\ref{G1def}) on $C(\rr)$:
$$
G_1 (\rr ) = \prod_{i = 0}^ k [(s_i - 1)! - 1 + \dd\sb
{1,s_i}];
$$
use this together with the definition of $w$ to get the inequality
$$
C(\rr ) / w (\rr ) \leq G_1 (\rr )
    \left [ \prod_{i=0}^ k Q(n_i , s_i )^ {-1} \right ].
$$
Let $Z_i$ be defined to be $1$ when $s_i = 1$ and $Q(n_i , s\sb
i)^ {-1}$
otherwise.  Then $\prod_i Q(n_i , s_i )^ {-1} \leq \prod_i
(Z_i (n_i + 3)
/ (n_i + 1)) \leq (\prod_i Z_i ) [(n+3)! / (2!(n+1)!)] = \frac{1}{2}
(n+3) (n+2)
\prod_i Z_i$.  So
\begin{equation} \label{quotient}
C (\rr ) / w(\rr ) \leq \frac{1}{2} (n+3)^ 2 D(\rr ),
\end{equation}
where $D(\rr) := G_1 (\rr ) \prod_{i=0}^ k Z_i$.

Now if there are $i<j$ for which $s_i \geq 5$ and $s_j \geq 4$,
consider the (lexicographically) first such pair $i,j$, and
let $\phi (\rr )$ be the sequence for which $s_i (\phi (\rr )) = 4$,
$s_j (\phi (\rr )) = s_j
+s_i - 4$, and $s_l (\phi (\rr)) = s_l$ for all $l \neq i,j$, where
$s_l$ always refers to $s_l (\rr )$ unless otherwise noted.  Then (ii)
of
Proposition~\ref{factorials} shows that $G_1 (\phi (\rr )) \geq G_1
(\rr )$,
with strict inequality if $s_j \neq 4$.
But also $\prod_i Z_i(\phi(\rr)) \geq \prod_i Z_i(\rr)$.  To see
this, note that
each factor $Z_l$ is increasing in $n_l$ (even when $s_l$ varies, as
long as
$s_l$ stays at least $4$ or else stays constant at $1$ or constant at
$3$).
Changing $\rr$ to $\phi (\rr )$ changes the $s_l$'s only by switching two
of
them that are both at least $4$, and furthermore the values of $n_l$ for
$\phi
(\rr )$ are equal to the values of $n_l$ for $\rr$ except when $i < l \leq
j$.
Since $n_l$ is increased for $i < l \leq j$, it follows that
$D(\phi (\rr ))
\geq D(\rr )$.

Let $\rr^ {(1)}$ be gotten from $\rr$ by iterating $\phi$ until
it is no longer possible to do so.  Note that $\rr^ {(1)}$ has
at most one $s_i$ greater than $4$ and that, by induction,
$D(\rr^ {(1)}) \geq D(\rr )$.

Now if $\max_i s_i (\rr^ {(1)} ) = 4$,
then let $\rr^ {(2)} = \rr^ {(1)}$; otherwise, if $s_i (\rr\sp
{(1)}) > 4$,
let $\rr^ {(2)}$ be the sequence for which $s_i (\rr^ {(2)} ) =
s_k (\rr^ {(1)} )$, $s_k (\rr^ {(2)} ) = s_i (\rr^ {(1)}
)$,
and $s_l (\rr^ {(2)} ) = s_l (\rr^ {(1)} )$ for $l \neq i,k$.
Clearly, $G_1 (\rr^ {(2)}) = G_1 (\rr^ {(1)})$.  For each $l
\neq i,k$,
$s_l$ is unchanged while $n_l$ is increased or remains the same,
and so $Z_l$ is not decreased.
Since $Z_i (\rr^ {(2)}) \geq (n+3)^ {-2} Z_i (\rr^ {(1)})$
and
$Z_k (\rr^ {(2)}) \geq Z_k (\rr^ {(1)})$, it follows that
$D(\rr^ {(2)}) \geq (n+3)^ {-2} D(\rr^ {(1)})$.  Summarizing the
progress so far, there is a sequence $\rr^ {(2)}$ with $s_i (\rr\sp
{(2)})
\leq 4$ for $i < k$ such that
$$
C(\rr ) / w (\rr ) \leq \frac{1}{2} (n+3)^ 2 D(\rr ) \leq \frac{1}{2}
(n+3)^ 2
D(\rr^ {(1)}) \leq \frac{1}{2} (n+3)^ 4 D(\rr^ {(2)}).
$$

The last step is to compare $D(\rr^ {(2)})$ to $D(\rr_0) = [(n-k)! -
1] (n-k+3)
(n-k+2)$, where $\rr_0$ is the sequence $(0 , 1 , 2 , \ldots , k , n+1)$.
Let $k_3 = k_3 (\rr^ {(2)})$ be the number of $i<k$ for which
$s_i (\rr^ {(2)}) = 3$ and let $k_4 = k_4 (\rr^ {(2)})$ be the
number of $i<k$ for which $s_i (\rr^ {(2)}) = 4$.  Then
$s_k (\rr^ {(2)}) = n+1 - k - 3 k_4  - 2 k_3$.  Now
the inequalities $(a!-1) / (b! - 1) > a!/b! > (a/e)^ {a-b}$ for $a>b>0$
from (iii) of Proposition~\ref{factorials} give
\begin{eqnarray*}
&&D(\rr_0) / D(\rr^ {(2)}) \\[2ex]
& \geq & [(n-k)! - 1] / [5^ {k_4} (n+3)^ {2k_4} (n+3)\sp
{k_3}
    ((n-k-3k_4-2k_3)!-1)] \\[2ex]
& \geq & ((n-k)/e)^ {3k_4 + 2k_3} / [5^ {k_4} (n+3)\sp
{2k_4 + k_3}]
 \\[2ex]
& \geq & (5 (n+3)^ {2/3})^ {3k_4 + 2k_3} / [5^ {k_4}
(n+3)^ {2k_4 +
 k_3}] \\[2ex]
& \geq & 5^ {2k_4 + 2 k_3} (n+3)^ {k_3 / 3} \\[2ex]
& \geq & 1,
\end{eqnarray*}
where the inequality $(n-k)/e \geq 5 (n+3)^ {2/3}$ comes from the
hypothesis of the lemma.  This proof of $D(\rr_0) \geq D(\rr^ {(2)})$
has
assumed $n - k - 3k_4 - 2k_3 \neq 1$, but can easily be modified to
treat the
contrary case.

Putting all this together gives
$C(\rr ) / w(\rr ) \leq n^ 6 (n-k)!$ for all $\rr$,
proving the $f(n,k)$ part of the lemma.

Finally, using Lemma~\ref{large k}, which is valid because $n - \lc 5 e
(n+3)^ {2/3} \rc + 1 \geq n - n^ {3/4} / 2$,
\begin{eqnarray*}
F(n,k)&\leq&\sum_{n - 5 e (n+3)^ {2/3} \geq l \geq k} f(n,l) + F(n,
n - \lc 5 e
(n+3)^ {2/3} \rc + 1) \\
&\leq&2 n^ 6 (n-k)! + \lc 5 e (n+3)^ {2/3} \rc (2n^ {7/8})\sp
{\lc 5 e
(n+3)^ {2/3} \rc - 1}.
\end{eqnarray*}
$\Cox$

In order to prove Lemma~\ref{Fsub1}, the following comparison
is needed between overlaps of pairs of paths connecting
$\0$ to $\1$ and overlaps of pairs of paths connecting
$x_i$ to $y_i$.

\begin{lem} \label{supplement}
Let $x_1, x_2$ be distinct vertices at level $L$ of $\B_{n+2L}$
and let $y_1, y_2$ be distinct vertices at level $n+L$ of $\B\sb
{n+2L}$,
with $x_1$ below $y_1$ and $x_2$ below
$y_2$.  Fix an oriented path $\gamma$ connecting $x_1$ to $y_1$.
For $0 \leq k \leq n$ and $i=1,2$, let $A_i(\gamma, k)$ be the set of
oriented paths $\gamma'$ from $x_i$ to $y_i$ that share exactly
$k$ edges with $\gamma$.  For $i=1,2$, let $A_i$ [$ = \cup_{k=0}^ n
A_i(\gamma,k)$] be the set of all oriented paths from $x_i$ to $y_i$.
Then there is a bijection $\phi$ from $A_1$ to $A_2$ such that
{\rm (i)} the set of edges that
$\phi (\gamma')$ has in common with $\gamma$ is a subset of the edges
that $\gamma'$ has in common with $\gamma$;
{\rm (ii)} hence if
$\gamma' \in A_1 (\gamma , k)$ then $\phi (\gamma') \in
A_2 (\gamma , j)$ for some $j \leq k$; and
{\rm (iii)} if
$\gamma'$ and $\gamma$ share either their first or last
edge, then the inclusion in (i) (and hence the inequality in (ii)) is strict.
\end{lem}

\noindent{\bf Proof:}  Viewing vertices of $\B_{n+2L}$ as subsets of
$\{ 1 , \ldots , n+2L \}$, the path $\gamma$ is represented by
a permutation of $y_1 \setminus x_1$.  Assume without loss of
generality that $x_1 = \{ 1 , \ldots , L\}$ and $y_1 = \{ 1 , \ldots ,
n+L \}$, and that the permutation representing $\gamma$ is in
fact the increasing permutation $(L+1 , \ldots , L+n)$.  We shall express the
desired $\phi$ as a bijection, call it $h$, between permutations
of $\{ L+1, \ldots, L+n \}$ and permutations of $y_2 \setminus x_2$.

Let $\gamma'$ connecting $x_1$ to $y_1$ be represented by a permutation
$\pi = (\pi (1) , \ldots , \pi (n))$ of $\{ L+1 , \ldots , L+n \}$.  We need to
define a corresponding permutation $h(\pi ) = (h(\pi)(i), i=1,\ldots,n)$
of $y_2 \setminus x_2$.  The idea is this:
$h$ tries to copy $\pi$, but is required only to copy
elements of $y_2 \setminus x_2$; so it replaces in corresponding
order the elements of $\{ L+1 , \ldots , L+n \}$ that are not elements
of $y_2 \setminus x_2$ by elements of $y_2 \setminus x_2$
that are not elements of $\{ L+1 , \ldots , L+n \}$.

The construction of $h$ can be expressed more formally as follows.  Let
$I :=
\{1, \ldots, n\} \cap [(y_2 \setminus x_2) - L]$ denote the set of
indices $i$,
$1 \leq i \leq n$, such that $L + i \in y_2 \setminus x_2$; thus the
points $L
+ i$, $i \in I$, form the intersection of $\{L + 1, \ldots, L + n\}$ and
$y_2
\setminus x_2$.  There are $m := n - |I|$ elements in $\{L + 1, \ldots,
L +
n\}$ that are not in $y_2 \setminus x_2$; label them in increasing
order as
$k_1 < k_2 < \cdots < k_m$.  Label the $m$ elements that are in
$y_2 \setminus
x_2$ but not in $\{L + 1, \ldots, L + n\}$ as $k_1' < k_2' <
\cdots < k_m'$.

Let $\pi = (\pi(1), \ldots, \pi(n))$ be a permutation of $\{L + 1, \ldots, L
+
n\}$.  If $\pi(i) \in y_2 \setminus x_2$, i.e., if $\pi(i) - L \in I$,
then let
$h(\pi)(i) = \pi(i) \in L + I$.  If $\pi(i) \notin y_2 \setminus x_2$,
i.e., if
$\pi(i) = k_t$ for some $1 \leq t \leq m$, then let $h(\pi)(i) = k_t'$.
It is easy to see that this yields a bijection between permutations
of $\{ L+1 , \ldots , L+n \}$ and permutations of $y_2 \setminus x_2$.
It
remains to show that it has the required properties.

Let $\gamma' \in A_1$ be represented by $\pi$ and $\phi(\gamma') \in A\sb
2$ by
$h(\pi)$.  Then $i$ belongs to $\rr(\gamma')$ (defined as a subsequence of
$(1,
\ldots, n)$ in the obvious fashion) if and only if
\begin{equation} \label{harder a}
\{ \pi(1), \ldots , \pi (i-1)\}
= \{ L+1, \ldots, L+i-1 \} \mbox{\ \ and\ \ }
\pi(i) = L + i.
\end{equation}
On the other hand, $i \in \rr (\phi (\gamma'))$
if and only if
\begin{equation} \label{harder b}
\{ h(\pi) (1) , \ldots , h(\pi ) (i-1) \}
= \{ 1, \ldots, L+i-1 \} \setminus x_2 \mbox{\ \ and\ \ }
h(\pi) (i) = L+i;
\end{equation}
in particular, $x_2 \subset \{1 \ldots, L+i-1\}$ is necessary
for~(\ref{harder
b}).

We claim that $\rr(\phi(\gamma')) \subseteq \rr(\gamma')$, i.e., that for each
$i \in \{1, \ldots, n\}$, (\ref{harder b}) implies~(\ref{harder a}).  Indeed,
suppose~(\ref{harder b}) holds for a given value of $i$.  Since
$h(\pi)(i) = L +
i \in \{L + 1, \ldots, L + n\}$, we have $\pi(i) = h(\pi)(i) = L + i$ by our
construction of $h(\pi)$.  Also, $\{h(\pi)(1), \ldots, h(\pi)(i-1)\}$ is the
union of some subset of $\{L + 1, \ldots L + i - 1\}$ and an initial
segment of
$\{k_1', \ldots, k_m'\}$ (namely, $\{1, \ldots, L\} \cap \{k_1',
\ldots,
k_m'\}$).  Therefore $\{\pi(1), \ldots, \pi(i-1)\}$ is the union of some
subset
of $\{L + 1, \ldots, L + i - 1\}$ and an initial segment of $\{k_1,
\ldots,
k_m\}$, and so equals $\{L + 1, \ldots, L + i - 1\}$.
Thus~(\ref{harder a}) is
established, and the proof of (i) is complete.  (ii) follows immediately from
(i).

To finish the proof of the lemma,
observe that since $x_1$ and $x_2$ are distinct, $\phi (\gamma')$
can never share its first edge with $\gamma$.  Similarly, $y_1$
and $y_2$ are distinct, so $\phi (\gamma')$ can never share its
last edge with $\gamma$.  Thus if $\gamma'$ shares either its
first or last edge with $\gamma$, the inclusion $\rr (\phi (\gamma'))
\subseteq \rr (\gamma')$ must be strict.   $\Cox$

\noindent{\bf Proof} of Lemma~\ref{Fsub1}:  Fix $x_i, y_i$ ($i=1,2$)
and
a path $\gamma$ connecting $x_1$ and $y_1$.
Use Lemma~\ref{supplement} to get a bijection $\phi$ from
paths connecting $x_1$ and $y_1$ to paths connecting $x_2$ and
$y_2$ with the properties stated therein.  Now the interval
in $\B_{n+2L}$ between $x_1$ and $y_1$ is isomorphic to $\B_n$.
Hence the number of paths $\gamma'$ connecting $x_1$ and $y_1$
and sharing at least $k$ edges with $\gamma$ is just $F(n,k)$, and
Lemmas~\ref{small k}~--~\ref{middle k} may be used to bound this.
Now since $\phi$ is a bijection and $\phi (\gamma')$ shares at
most as many edges with $\gamma$ as $\gamma'$ does, this immediately
gives $H (n,L,k,x_1,x_2,y_1,y_2) \leq F(n,k)$; maximizing over
$x_1$, $y_1$, $x_2$, and $y_2$ and applying Lemmas~\ref{large k}
and~\ref{middle k} establishes~(\ref{large
Fsub1})~and~(\ref{middle Fsub1}).
To show~(\ref{small Fsub1}), note that if $\gamma''$ connects
$x_2$ to $y_2$ and shares at least $k$ edges with $\gamma$ then
$\gamma' = \phi^ {-1} (\gamma'')$ either shares strictly more
edges with $\gamma$ or shares exactly the same edges in which
case the shared edges include neither the first nor last edge.
The number of $\gamma'$ in the former category is at most
$F(n,k+1)$, while the number of $\gamma'$ in the latter category
is at most $o((k+1)(n-k)!)$ uniformly in $0 < k \leq K(n)$ according
to the last part of Lemma~\ref{small k}.  But $F(n,k+1) \leq
(1+o(1)) (k+2) (n-k-1)! = o((k+1)(n-k)!)$ uniformly in $0 < k \leq K(n)$
according to the first part of Lemma~\ref{small k},
and the desired conclusion follows.    $\Cox$

\section{Oriented percolation and oriented first-passage percolation}

\subsection{Oriented percolation}

The first application of Lemmas~\ref{2nd moment method}~--~\ref{Fsub1}
will be to ordinary oriented percolation.  This means that the edges
are independently open with probability $p$ and closed otherwise.
The problem considered here is to determine how large $p$ should be as
a function of the dimension $n$ in order that $\0$ and $\1$ may
be connected with substantial probability.
It turns out that the answer is $p = e/n$
in the sense that if $p = c/n$ then for $c < e$ the probability
that $\0$ is connected to $\1$ tends to $0$ as $n \rightarrow
\infty$, while for $c \geq e$, the probability
that $\0$ is connected to $\1$ approaches a positive limit which we calculate.
(It is easy to see that the limit is not $1$ since the disconnection
probability is at least the chance that $\0$ is isolated, namely, $(1 -
c/n)^ n \to e^ {-c}$.)  In particular, the limiting connection
probability is
$\doteq 0.8416$ at the critical value $c = e$.  We consider oriented
percolation
mainly as a warm-up to the arguments used to analyze oriented first-passage
percolation.  We do not intend for our results to be viewed as a
complete analysis of the threshold behavior of oriented percolation with edge
probabilities $c/n$ as $c$ crosses the critical value $e$.

The next lemma derives a lower bound of $(e-1)^ 2/e^ 2 \doteq 0.400$
on the
percolation probability in the critical region by using the unenhanced second
moment method.  The bound is not tight but will be a crucial ingredient to the
proof of the sharper result.

\begin{lem} \label{e case}
Fix $L \in \R$.
Let each edge of $\B_n$ be independently open with probability
$p = e/(n+2L)$.
Then
$$
\liminf_n \P(\0 \mbox{\rm \ is connected to $\1$ by an oriented open
path}) \geq
\frac{(e-1)^ 2}{e^ 2}.
$$
\end{lem}

\noindent{\bf Proof:}  Observe that $\E N^ 2$ is the
sum over pairs of paths of the probability that both paths
are open.  By symmetry this is the same as $n!$ times the
sum when the first path is fixed, say as the path $\gamma$
corresponding to the identity permutation.  The probability
that both $\gamma$ and $\gamma'$ are open depends only on
the number $k$ of edges they share and is equal to $p^ {2n-k}$.
Thus by Lemma~\ref{2nd moment method} it suffices to find a
finite upper bound for
\begin{equation} \label{moment estimate}
{\E N^ 2 \over (\E N)^ 2} = n! {\sum_k f(n,k) p^ {2n-k} \over
    (p^ n n!)^ 2 } = \sum_k [f(n,k) (\frac{n+2L}{e})^ k / n!].
\end{equation}

We bound this in three pieces, corresponding to
Lemmas~\ref{small k}~--~\ref{middle k}; of course we may assume $n$ to be as
large as needed.  When $k < 12 \ln n$, the summand in the final sum
of~(\ref{moment estimate}) is by Lemma~\ref{small k} at most $(1+o(1))
(k+1)
e^ {-k} (n+2L)^ k / [n(n-1) \cdots (n-k+1)]$ uniformly in $k$ as
$n \rightarrow
\infty$. This is at most $(1+o(1)) (k+1) [(n+2L)/(n-12 \ln n)]^ k
e^ {-k}$
and so the sum over $k$ is at most $(1+o(1)) \left( e / \left[ e
- (n+2L)/(n-12\ln n) \right] \right)^ 2 = (1 + o(1)) e^ 2 / (e-1)\sp
2$.

The contribution from any term with $k > n - n^ {3/4} / 2$ is
at most
$$
(n-k+1) (2n^ {7/8})^ {n-k} [(n+2L)/e]^ k / n!
$$
according to Lemma~\ref{large k}.  Using the inequality $n! > (n/e)^ n
n^ {1/2}$
gives an upper bound for these terms of $(n-k+1) (2e n^ {-1/8}
\frac{n}{n+2L})^ {n-k} n^ {-1/2} (1 + \frac{2L}{n})^ n$.  The sum
over $k$ is then
at most $[1/(1-2en^ {-1/8} \frac{n}{n+2L})]^ 2 n^ {-1/2}
(1 + \frac{2L}{n})^ n =
(1 + o(1)) e^ {2L} n^ {-1/2} = o(1)$ as $n \to \infty$.

Finally, consider the contribution from terms with $12 \ln n
\leq k \leq n-5e(n+3)^ {2/3}$.  By Lemma~\ref{middle k} this is at most
$$
\sum_{12 \ln n \leq k \leq n-5e(n+3)^ {2/3}} g(n,k),
$$
where $g(n,k) = n^ 6 (n-k)! [(n+2L)/e]^ k / n!$.   Now
$g(n,k+1) / g(n,k) = (n+2L) / (e(n-k))$, which is increasing in $k$;
hence $g$ is  U-shaped in $k$ for fixed $n$ (and $L$), with its minimum
at $1 + \lfloor n - (n+2L)/e \rfloor$.  In particular, the maximum
of $g(n,k)$ for fixed $n$ over an interval of values of $k$
is achieved at an endpoint.  Thus $\max \{ g(n,k) : 12 \ln n
\leq k \leq n - 5e(n+3)^ {2/3} \}$ is achieved at an endpoint.  At the
first endpoint, $k = \lc 12 \ln n \rc$ and $g(n,k) =
n^ 6 e^ {-\lc 12 \ln n \rc} \prod_{i=0}^ {k-1}
[(n+2L)/(n-i)] \leq n^ 6 n^ {-12}
[(n+2L)/(n-12 \ln n)]^ {\lc 12 \ln n \rc} = (1 + o(1)) n^ {-6}$.
At the second
endpoint, $k = n - \lc 5e(n+3)^ {2/3} \rc$ and the inequality $n! >
(n/e)^ n$
gives
\begin{eqnarray*}
g(n,k)&  = & n^ 6 \lc 5e(n+3)^ {2/3} \rc !\,[(n+2L)/e]\sp
{n - \lc 5e(n+3)^ {2/3}
              \rc} / n! \\
      &  < & n^ 6 \lc 5e(n+3)^ {2/3} \rc !\,[(n+2L)/e]\sp
{-\lc 5e(n+3)^ {2/3} \rc}
              (1 + 2Ln^ {-1})^ {n - \lc 5e(n+3)^ {2/3} \rc} \\
      &  = & (1+o(1)) e^ {2L} n^ 6 (\lc 5e(n+3)^ {2/3}
\rc / e)^ {\lc 5e(n+3)^ {2/3}
              \rc} \sqrt{2\pi 5en^ {2/3}} [e/(n+2L)]\sp
{\lc 5e(n+3)^ {2/3} \rc} \\
      &  = & (1+o(1)) e^ {2L} n^ 6 \left( \frac{\lc 5e(n+3)\sp
{2/3} \rc}{n+2L}
              \right)^ {\lc 5e(n+3)^ {2/3} \rc} \sqrt{2\pi
5en^ {2/3}}
\end{eqnarray*}
by Stirling's formula.
This clearly tends to $0$ faster than any power of $n$ and is thus
$o(n^ {-6})$.
Now the terms at both ends have been shown to be at most $(1 + o(1))
n^ {-6}$ and
there are at most $n$ terms, so the total contribution from these terms is at
most $(1 + o(1)) n^ {-5} = o(1)$.

To sum up, for $n$ sufficiently large, $\sum_k f(n,k) [(n+2L)/e]\sp
k / n!$ is at
most the sum of the contributions from the three ranges,
which was computed to be $(1+o(1)) e^ 2 / (e-1)^ 2 + o(1) + o(1)$.
Thus by Lemma~\ref{2nd moment method}, $\P(\0
\mbox{ is connected to } \1) \geq (1+o(1)) (e-1)^ 2 / e^ 2$.  $\Cox$

\begin{th} \label{ordinary}
Let each edge of $\B_n$ be independently open with probability
$p=c/n$.  Then
$\P(\0 \mbox{\rm \ is connected to $\1$ by an oriented open path})$
converges to a limit as $n \rightarrow \infty$.
The limit is $0$ if $c < e$ and is $(1-x(c))^ 2$ if $c \geq e$, where
$x(c)$ is
the extinction probability for a Poisson($c$) Galton--Watson process, namely,
the solution in $(0,1)$ to $x=e^ {c(x-1)}$.
\end{th}

Note that as $c \to \infty$, $x(c) = (1+o(1)) e^ {-c} = o(1)$, so that the
limiting connection probability is $1 - (1+o(1)) 2e^ {-c} \to 1$.

\noindent{\bf Proof:}  There are $n!$ oriented paths from $\0$ to $\1$.
Let $N$ be the random number of these that consist entirely of
open edges.  For each path $\gamma$ the probability that $\gamma$
is open is $p^ n$, so $\E N = n! p^ n$.  If $c < e$ then
$\E N = n! c^ n n^ {-n} = (1+o(1)) (c/e)^ n \sqrt{2 \pi n}$, which
tends to $0$.  Since $\P (N > 0) \leq \E N$, this proves the
first part.

For the second part, fix $c \geq e$.  Also fix $\ee > 0$.
Write $M = \lceil 1 / \ee \rceil$.
For $i = 1,2,\ldots$, let $A_i$ be the set of vertices at level $i$
reachable from $\0$ in $\B_n$.
For any fixed $i$, as $n \rightarrow \infty$, the
joint distribution of $|A_0| , \ldots, |A_i|$ approaches in total
variation the distribution of a Galton--Watson process with
the number of offspring of each particle Poisson distributed with
mean $c$.  Because a surviving branching process proliferates, an integer
$L =
L(\ee)$ may be chosen so that $\P(|A_L| \geq M) \geq (1-x(c)) - \ee$ for
sufficiently large $n$, where $x(c)$ is the extinction
probability for the Galton--Watson process, namely, the solution
in $(0,1)$ to $x=e^ {c(x-1)}$.  Let $B_j$ be $A_j$ upside down,
i.e., the set of vertices at distance $j$ from $\1$ that can reach
$\1$.  Then by symmetry and independence, we have $\P(F) \geq (1-x(c))\sp
2
- 2\ee$, where $F$ is the event $\{|A_L| \geq M \mbox{ and } |B_L|
\geq M\}$.
Now if either of the two sets $A_L$ or $B_L$ is
empty, then $\0$ is not connected to $\1$, so the $\limsup$ of $(1 - x(c))\sp
2$ is
established by the convergence in total variation.

For the lower bound we employ the enhanced second moment method described
following Lemma~\ref{2nd moment method}, although the details here are
slightly
different.  The variance-absorbing random variable $Z$ is $\min(|A_L|,
|B_L|)$.  Uniformly in $z \geq M$, we show
\begin{equation} \label{N given Z}
\P(N > 0 | Z = z) \geq (1 - o(1)) \left[ 1 + \frac{e^ 2}{M (e-1)^ 2}
\right]^ {-1}
\end{equation}
as $n \to \infty$.  Then
$$
\P(N > 0) \geq (1 - o(1)) \left[ 1 + \frac{e^ 2}{M (e-1)^ 2}
\right]^ {-1} \P(F)
$$
and so
$$
\liminf_n \P(N > 0) \geq \left[ 1 + \frac{e^ 2}{M (e-1)^ 2}
\right]^ {-1} [(1 -
x(c))^ 2 - 2 \ee].
$$
Letting $\ee \downarrow 0$ gives $\liminf_n \P(N > 0) \geq (1 - x(c))\sp
2$, as
desired.

Henceforth tacitly conditioning on $Z = z \geq M$, we prove~(\ref{N given Z})
by applying the second moment method to a truncation $N'$ of $N$, as follows.
First, reduce the probability that any given edge falling between levels $L$
and $n - L$ is open from $p = c/n$ to $p = e/n$.  The obvious coupling
argument
shows that this diminishes $N$ stochastically.
Then let $x_1 , \ldots , x_M$ be an enumeration of the first $M$
vertices
of $A_L$ and let $y_1 , \ldots , y_M$ be an enumeration of
the first $M$ vertices of $B_L$, in some arbitrary ordering of the
vertices at
levels $L$ and $n - L$, respectively.  Let $N_i$ be the number
of open paths connecting $x_i$ to $y_i$, and consider $N' := \sum\sb
{i=1}^ M
N_i$.  Essentially the same
calculations that showed $\E N^ 2 / (\E N)^ 2$ to be bounded in the
proof of
Lemma~\ref{e case} will be used to show $\E (N')^ 2 / (\E N')^ 2
\leq (1+o(1))
\left( 1 + e^ 2 / [M (e - 1)^ 2] \right)$ via the two results
$\E N_i^ 2 / (\E N_i)^ 2 \leq (1 + o(1)) e^ 2/(e-1)\sp
2$, uniformly in $i$, and
$\E(N_i N_j) / [(\E N_i) (\E N_j)] \leq (1 + o(1))$, uniformly
in pairs $i \neq
j$.  Indeed, these last two statements imply $\E(N')^ 2 = \sum_i \E
N_i^ 2 +
\sum_{i \neq j} \E(N_i N_j)
\leq (1 + o(1)) [e^ 2 / (e-1)^ 2] \sum_i (\E N_i)^ 2
+ (1+o(1)) \sum_{i \neq j} (\E N_i) (\E N_j)
\leq (1 + o(1)) M [e^ 2 / (e-1)^ 2] (\E N_1)^ 2 + (1 + o(1)) M
(M - 1) (\E N_1)^ 2
\leq (1+o(1)) \left( 1 + e^ 2 / [M (e - 1)^ 2] \right) (\E N')\sp
2$.  Now simply
apply Lemma~\ref{2nd moment method} to deduce~(\ref{N given Z}) for $N'$ and
hence for $N$.

Working for convenience with $p = e / (n + 2L)$ in $\B_{n + 2L}$ rather
than
with $p = e/n$ in $\B_n$, and given distinct vertices $x_1$ and $x_2$
in level
$L$ and distinct vertices $y_1$ and $y_2$ in level $n + L$, let $N_i$
denote
the number of open paths from $x_i$ to $y_i$, $i=1,2$.  We must show $\E
N_1^ 2
/ (\E N_1)^ 2 \leq (1 + o(1)) e^ 2 / (e-1)^ 2$ and $\E(N_1 N\sb
2) / [(\E N_1) (\E
N_2)] \leq 1 + o(1)$ as $n \to \infty$, uniformly in the choice of $x_1,
x_2,
y_1, y_2$.  Now the interval
from $x_1$ to $y_1$ is isomorphic to $\B_n$, so $N_1$ has the same
distribution
as the total number of open paths, $N$, in the proof of Lemma~\ref{e case}.
Thus we immediately obtain $\E N_1^ 2 / (\E N_1^ 2) \leq (1 +
o(1)) e^ 2 / (e -
1)^ 2$ (uniformly in $x_1$).  For the other inequality, mimic the
calculation
from Lemma~\ref{e case} to get
$$
{\E(N_1 N_2) \over (\E N_1) (\E N_2)} \leq n! {\sum_k
H(n,L,k,x_1,x_2,y_1,y_2)
\left( \frac{e}{n+2L} \right)^ {2n-k} \over \left[ \left(
\frac{e}{n+2L}
\right)^ n n! \right]^ 2} \leq \sum_k F_1(n,L,k) \left(
\frac{n+2L}{e} \right)^ k /
n!.
$$
Now break the sum into three pieces again, corresponding to values
$k < 12 \ln n$, $k > n - n^ {3/4}/2$, and all $k$ in between.
This time use Lemma~\ref{Fsub1} instead of Lemmas~\ref{small
k}~--~\ref{middle k}, one difference being that the contribution
for terms with $0 < k < 12 \ln n$ is now $o((k+1)(n-k)!) [(n+2L)/e]^ k /
n!$,
and thus the contribution from those terms and the unit contribution from
the $k=0$ term sum to $1+o(1)$.  As before, the contribution from large
$k$ is $o(1)$, and the contribution from the first term in the
bound~(\ref{middle Fsub1}) for intermediate values of $k$ is also $o(1)$.
We
finish the proof of the theorem by showing that the contribution from the
second
term in~(\ref{middle Fsub1}) is also $o(1)$:
\begin{eqnarray*}
\lefteqn{\sum_{12 \ln n \leq k \leq n - 5 e (n+3)^ {2/3}} \lc 5 e
(n+3)^ {2/3}
\rc (2 n^ {7/8})^ {\lc 5 e (n+3)^ {2/3} \rc - 1} [(n+2L)/e]^ k /
n!}  \\
 &\leq& n \lc 5 e (n+3)^ {2/3} \rc (2 n^ {7/8})^ {\lc 5 e (n+3)\sp
{2/3} \rc - 1}
[(n+2L)/e]^ {n - \lc 5 e (n+3)^ {2/3} \rc} / n! \\
 &  = & (1+o(1)) \frac{5 e^ {2L+1}}{2 \sqrt{2 \pi}} n^ {7/24}
\left(
\frac{2e}{n^ {1/8}} \right)^ {\lc 5 e (n+3)^ {2/3} \rc},
\end{eqnarray*}
which vanishes at a rate faster than any power of $n$.
       $\Cox$

\subsection{Oriented first-passage percolation}

Now consider oriented first-passage percolation (OFPP).
Give each edge in $\B_n$ an upward orientation and assign $\iid$
random variables $X_e$ with common density $f$ to each edge $e$.
The problem in OFPP is to determine the minimum value over
oriented paths from $\0$ to $\1$ of the sum along the path of
the $X_e$'s.  Under mild conditions on $f$, it turns out
(Theorem~\ref{OFPP time}) that this random minimum converges in probability
to
$1 / f(0)$ as $n \rightarrow \infty$.  By multiplying
every edge-passage time $X_e$ by a constant, it can be assumed without loss
of generality (provided $0 < f(0) < \infty$) that $f(0) = 1$.

As Aldous (1989) points out, use of the
exponential distribution $f(x) = e^ {-x}$ simplifies some of the
calculations
involved but is not necessary.  The following two lemmas, treating the
exponential and more general cases, respectively, produce large deviation
estimates that correspond to the probability $[e/(n+2L)]^ {2n-k}$ in the
proof of
Theorem~\ref{ordinary}.

%*% The proof of Lemma~\ref{deviation} has been revised.
\begin{lem} \label{deviation}
\item[{\rm (i)}] Let $S_n$ be the sum of $n \geq 1$ $\iid$
random variables $Y_i$, each exponential with mean $1$, and let
$u$ be a real number in $[0,1]$.  Then $\P (S_n \leq u) =
(1 + K_1(u,n)) e^ {-u} u^ n / n!$ with $0 \leq K_1(u,n)
\leq e / (n+1) \leq 2$.
 \item[{\rm (ii)}] Given $1 \leq k \leq n-1$, let $S_n' = \sum_{i=1}\sp
k Y_i +
\sum_{i=k+1}^ n Y_i'$, where $Y_1, \ldots, Y_n, Y_{k+1}',
\ldots, Y_n'$ are
independent and identically distributed.  Then $\P (S_n \leq
1\mbox{\rm\ and\
}S_n' \leq 1) \leq K_2 R(n,k)$, where
$$
R(n,k) =  2^ {2n-2k} e^ {2n-k} (2n-k)^ {-(2n-k)} /
    \left[ (n-k)^ {1/2} (2n-k)^ {1/2} \right]
$$
and $K_2$ is constant.  Furthermore, for $1 < k \leq n - 1$, $R(n,k-1) /
R(n,k)
\leq K_3 / n$ for some constant  $K_3$.
\end{lem}

\noindent{\bf Proof:}  The key for (i) is the standard switching relation
$\P(S_n \leq u) = \P(X_u \geq n)$, where $X = (X_u)_{u \geq 0}$
is a Poisson
process with unit intensity parameter.  Thus $\P(S_n \leq u) \geq \P(X\sb
u = n) =
e^ {-u} u^ n / n!$.  Moreover,
\begin{eqnarray*}
\P(S_n \leq u)
 &  = &\sum_{m=n}^ {\infty} e^ {-u} \frac{u^ m}{m!} = \frac{u\sp
n}{n!}
        \times e^ {-u} \sum_{l=0}^ {\infty} \frac{u^ l}{(n+l)
\cdots (n+1)} \\
 &\leq&\frac{u^ n}{n!} \times e^ {-u} \sum_{l=0}^ {\infty}
\frac{u^ l}{l!} =
        \frac{u^ n}{n!}.
\end{eqnarray*}
Using the switching relation together with this crude upper bound, we obtain
\begin{eqnarray*}
\P(S_n \leq u)
 &  = &\P(X_u = n) + \P(X_u \geq n+1) \\
 &  = &e^ {-u} \frac{u^ n}{n!} + \P(S_{n+1} \leq u) \\
 &\leq&e^ {-u} \frac{u^ n}{n!} (1 + e^ u \frac{u}{n+1}) \\
 &\leq&e^ {-u} \frac{u^ n}{n!} (1 + \frac{e}{n+1}),
\end{eqnarray*}
as desired.

For (ii), we both prove the large deviations inequality and show that it is
tight.  Begin by writing
$\P(S_n, S_n'\leq 1) = \int_0^ 1\,\P(S_k \in du)\,[\P(S\sb
{n-k} \leq 1 - u)]^ 2$.
Using the bounds from (i) gives $$
\P(S_n, S_n'\leq 1) = \int_0^ 1\,[e^ {-u} u^ {k-1} /
(k-1)!]
[e^ {2u-2} (1-u)^ {2n-2k}] [1 + K_1(1-u,n-k)]^ 2 / [(n-k)!]\sp
2\,du,
$$
and using lower and upper bounds for $e^ {u-2}$ and $K_1(1-u,n-k)$ bounds
this below by
$$
{e^ {-2} \over (k-1)! [(n-k)!]^ 2} \int_0^ 1\,u^ {k-1}
(1-u)^ {2n-2k}\,du
$$
and above by
$$
{9e^ {-1} \over (k-1)! [(n-k)!]^ 2} \int_0^ 1\,u^ {k-1}
(1-u)^ {2n-2k}\,du.
$$

Now the integral is equal to $(k-1)! (2n-2k)! / (2n - k)!$, so
$\P(S_n, S_n' \leq 1)$ is bounded between $e^ {-2}$ and $9 e\sp
{-1}$ times
${{2n - 2k} \choose {n-k}} / (2n-k)!.$
This can be approximated using Stirling's formula, for which it will suffice
to note that the error factor of $e^ {1/(12n)}$ is bounded.  Thus
$\P(S_n,
S_n' \leq 1)$ is bounded between positive constant multiples of $R(n,k)$.

To see that $R(n,k-1) / R(n,k) \leq K_3 / n$, note that the exact
quotient is $4e [(2n - k) / (2n - k +1)]^ {2n-k} (2n - k + 1)^ {-1}
[(n-k) / (n - k + 1)]^ {1/2} [(2n-k) / (2n - k + 1)]^ {1/2}$.
The factor $[(n-k) / (n - k + 1)]^ {1/2}$ is between $\sqrt{1/2}$ and $1$,
while the rest of the product equals $(1+o(1))\,4/(2n-k) \leq
(1+o(1))\,4/n$, uniformly in $k$.  This proves the claim and finishes that of
the lemma.  $\Cox$

The Lipschitz condition in the following lemma does not give the most
general $f$ for which the first-passage times can be calculated, but it does
cover most non-pathological cases.

\begin{lem} \label{nonexponential}
Let $f$ be a probability density on $[0,\infty)$,
and suppose that $f(0) = 1$ and that $f$ satisfies a ``global'' Lipschitz
condition at the origin: $|f(x) - f(0)| \leq K_4 x$ for some positive
$K_4 <
\infty$ and all $x \geq 0$. Let $T_n$ be the sum of $n$ independent random
variables $Z_1 , \ldots , Z_n$ with common density $f$ and let $S_n$
be the sum
of $n$ i.i.d.\ exponentials $Y_1 , \ldots , Y_n$ with unit mean.  Then
for
$0 < u \leq 1$, $\P (T_n \leq u) \leq e^ {(1+K_4) u}\,\P(S_n \leq
u)$
and $\liminf_n \P(T_n \leq u) / \P (S_n \leq u) \geq e^ {-K_4
u}$.
Similarly, if $T_n' = \sum_{i=1}^ k Z_i + \sum_{i=k+1}^ n
Z_i'$,
where $Z_1, \ldots, Z_n, Z_{k+1}', \ldots, Z_n'$ are i.i.d., then
$\P(T_n \leq
1\mbox{\rm\ and\ }T_n' \leq 1) \leq K_5 R(n,k)$, where $K_5 :=
e^ {2+2K_4} K_2$.
\end{lem}

\noindent{\bf Proof:}  For the upper bound, note that the Radon--Nikodym
derivative of $Z_i$ with respect to $Y_i$ at $x$ is $f(x) / e\sp
{-x}
\leq (1 + K_4 x) / e^ {-x} \leq e^ {(1+K_4)x}$.  Thus the
Radon--Nikodym
derivative of the $n$-tuple $(Z_1, \ldots, Z_n)$ with respect to
$(Y_1, \ldots, Y_n)$ at $(x_1, \ldots, x_n)$ is at most
$e^ {(1+K_4) \sum x_i}$, and hence the
derivative of $T_n$ with respect to $S_n$ at $x$ is at most
$e^ {(1+K_4)x}$.
This establishes the upper bound for $\P(T_n \leq u)$.  Together with
(ii) of
Lemma~\ref{deviation}, this argument also establishes the
upper bound for $\P (T_n \leq 1, T_n' \leq 1)$.

For the lower bound, we first establish the fact that for any fixed $u \in
(0,1]$ and $\dd \in (0,u)$, $\P(\max_i Y_i < \dd \| S_n \leq u)$
converges to
$1$ as $n \rightarrow \infty$.  To see this, note that $\P (S_n
\leq u , \max_i Y_i \geq \dd) \leq n \P(S_n \leq u , Y_1 \geq
\dd)
\leq n \P (S_{n-1} \leq u - \dd)$; part (i) of the
last lemma shows that for $n \geq 2$ this is at most $3 n e^ {\dd} e\sp
{-u}
(u-\dd)^ {n-1} / (n-1)!$.  Then the lower bound from (i)
of the previous lemma shows that $\P(\max_i Y_i \geq \dd \| S_n
\leq u)
\leq 3 e^ {\dd} ((u - \dd)/u)^ {n-1} n^ 2/u$ which vanishes at an
exponential rate
as $n \rightarrow \infty$, proving the claim.

Now, given $K_4' > K_4$, a lower bound for the Radon--Nikodym derivative
of
$Z_i$ with respect to $Y_i$ is $1 - K_4 x \geq e^ {-K_4' x}
\one_{\{x < \dd\}}$
for an appropriate $\dd \in (0,u)$.  Then $\P (T_n \leq u)
\geq e^ {-K_4' u} \P (S_n \leq u , \max_i Y_i < \dd)$; thus
$\P (T_n \leq u) / \P (S_n \leq u)$ is at least $e^ {-K_4' u} \P
( \max_i Y_i
< \dd \| S_n \leq u)$ and so has a $\liminf$ of at least $e^ {-K_4'
u}$ by
the fact in the previous paragraph.  Now let $K_4' \downarrow K_4$.
$\Cox$

We are now ready for the main theorem for OFPP.
\begin{th} \label{OFPP time}
Let the edges of $\B_n$ be assigned i.i.d.\ positive random
passage times with common density $f$, and assume that
$|f(x) - 1| \leq K_4 x$ for all $x \geq 0$.  Then
the oriented first-passage percolation time $T = T^ {(n)}$ for $\B_n$
converges
to $1$ in probability as $n \rightarrow \infty$.
\end{th}

\noindent{\bf Proof:}  Let $\ee$ be small and positive.  With
$X_{vw}$ being i.i.d.\ with common density $f$ and $\gamma$ an oriented path
from $\0$ to $\1$ in $\B_n$, let $T_n(\gamma )$ be the sum of
$X_{vw}$ along
edges $\vw$ of $\gamma$, so that the first-passage time $T$ is just
$\min_\gamma
T_n(\gamma )$.
Let $T_n$ denote a random variable distributed identically
to each $T_n (\gamma)$ and let $S_n$ denote the sum of
$n$ i.i.d.\ exponentials of unit mean, as in the lemmas.
Showing that $\P (T \leq 1-\ee ) \rightarrow 0$ is easy.  Let
$N$ be the number of $\gamma$ for which $T_n(\gamma ) \leq 1-\ee$.
Then, using Lemmas~\ref{nonexponential} and~\ref{deviation},
$\P(N > 0) \leq \E N = n! \P(T_n \leq 1 - \ee )
\leq n! \exp[(1 + K_4) (1 - \ee)] \P(S_n \leq 1 - \ee) \leq 3
\exp[K_4 (1 -
\ee)] (1 - \ee)^ n = o(1)$ as $n \to \infty$, where $S_n$ is the
sum of $n$
i.i.d.\ exponentials of mean $1$.

To show that $\P (T \leq 1+\ee )$ is bounded away from $0$, one
can mimic the proof of Lemma~\ref{e case}, but in order to
show that this probability converges to $1$, we need to find
another auxiliary random variable to reduce the variance.
It will be easier to work in $\B_{n+2}$.
Let $A_0$ be the random set of neighbors $v$ of $\0$ for which
the edge $\overline {\0 v}$ has $X_{\0 v} \leq \ee/2$.  Similarly,
let $A_1$ be the random set of neighbors $v$ of $\1$ for which
the edge $\overline {v \1}$ has $X_{v \1} \leq \ee/2$.  Let $b'$
be the minimum of $|A_0|$ and $|A_1|$.  Enumerate the elements
of $A_0$ by $x_1 , x_2 , \ldots$ and the elements of $A_1$ by
$y_1, y_2, \ldots$ in such a way that for $1 \leq i \leq b'$,
$y_i$ lies above $x_i$.  This is easy to do since there is
only one neighbor of $\1$ that does not lie above any given $x_i$.
Let $b = \lceil \ee n / 4 \rceil$.
The first thing to observe is that $b' > b$ with
probability converging to $1$ as $n \rightarrow \infty$.  This
is immediate from the fact that $b'$ is the minimum of two
independent random variables that are binomial with parameters $n+2$
and $\int_0^ {\ee/2}\,f(x)\,dx \geq \int_0^ {\ee/2}\,(1-K_4
x)\,dx > \ee/4$.  Now
condition on the event that $b' > b$.  It suffices to show that the
probability
of finding an oriented path $\gamma$ connecting $x_i$ to $y_i$ with
$T_n (\gamma
) \leq 1$ for some $i \leq b$ converges (with appropriate uniformity) to
$1$,
conditionally given $b' > b$ and the enumeration of the $x_i$'s and $y\sb
i$'s.
What will in fact be shown is that, uniformly over all choices of vertices
$x_1, x_2, \ldots, x_b$ neighboring $\0$ and $y_1, y_2, \ldots,
y_b$
neighboring $\1$ with $y_i$ above $x_i$ for each $i$, the probability of
finding
an oriented path of passage time at most $1$ connecting some $x_i$ to $y\sb
i$
tends to~$1$.

For this we use the second moment method.  For $1 \leq i \leq b$,
let $N_i$ be the number of paths connecting $x_i$ to $y_i$
with passage time at most $1$.  Let $N = \sum_{i=1}^ b N_i$.
The interval in $\B_{n+2}$ from $x_i$ to $y_i$ is isomorphic
to $\B_n$.  It is therefore easy, using Lemmas~\ref{deviation}
and~\ref{nonexponential}, to see that
$$
\E N = \sum_{i=1}^ b \E N_i = b\,\E N_1 = b\,n!\,\P(T_n \leq
1) = b\,c_0,
$$
where $c_0 = c_0 (n)$ is bounded between positive constants.
Now $\E N^ 2 = \sum_i \E N_i^ 2 + \sum_{i \neq j} \E (N_i
N_j)$.
If we can show that
\begin{equation} \label{twobounds}
\E N_1^ 2 = O(1)\mbox{\ \ \ and\ \ \ }\E(N_1 N_2) \leq (1 + o(1))
c_0^ 2,
\end{equation}
uniformly in the choice of $x_1, x_2, y_1, y_2$, then, it will
follow that
$\E N^ 2 / (\E N)^ 2 \leq O(b^ {-1}) + 1 + o(1)$, uniformly in the
choice of the
$2b$ vertices.  Since $b$ tends to infinity with $n$, this bound converges to
$1$, and so $\P(T \leq 1 + \ee) \to 1$, proving the theorem.

Each part of~(\ref{twobounds}) is established in pieces, in a manner
similar to the bounding of~(\ref{moment estimate}). For any fixed
$\gamma$ connecting $x_1$ to $y_1$, $\E N_1^ 2$ is given by $n!$
times the sum
over $\gamma'$ connecting $x_1$ to $y_1$ of $\P(T_n(\gamma) \leq 1
\mbox{\ and\
} T_n(\gamma') \leq 1)$. Break the sum into three ranges according to the
number
$k$ of edges shared by $\gamma$ and $\gamma'$ as before, and additionally
separate the cases $k=0$ and $k=n$.  The case $k=0$ means $T_n(\gamma )$ is
independent of $T_n (\gamma')$, so the contribution to $\E N_1^ 2$in
this case is
at most $(\E N_1)^ 2$; and the case $k=n$ has $\gamma = \gamma'$, so the
contribution in this case is exactly $\E N_1$. Using
Lemma~\ref{nonexponential}
for the other three ranges and recalling that the case $k=n-1$ is impossible,
the sum can be bounded by
\begin{eqnarray*}
\E N_1^ 2 &  \leq & (\E N_1)^ 2 + \E N_1 + n! \sum_{k=1}\sp
{n-2} f(n,k) K_5 R(n,k)
\\[2ex]
& \leq & (c_0^ 2 + c_0) \\
& & +\,(1+o(1)) n! \sum_{0 < k < 12 \ln n} (k+1) (n-k)! K_5 R(n,k) \\
& & +\,n! \sum_{n - 2 \geq k > n- n^ {3/4} /2} (n-k+1) (2n\sp
{7/8})^ {n-k}
    K_5  R(n,k) \\
& & +\,n! \sum_{12 \ln n \leq k \leq n - 5e(n+3)^ {2/3}} n^ 6
(n-k)!
    K_5  R(n,k)  ,
\end{eqnarray*}
for large enough $n$. Now it will not be too hard to show that this is
$c_0^ 2 +
c_0 + O(1) = O(1)$, but before doing so, notice how similar the above bound
is to
a good bound on $\E(N_1 N_2)$.  The $k=0$ term
for $\E(N_1 N_2)$ has the same bound as above, but the $k=n$ term
vanishes.  The necessary changes are completed by using $F_1$
in place of $f$.  This allows the $(k+1)(n-k)!$
to be replaced by $o((k+1)(n-k)!)$ according to Lemma~\ref{Fsub1},
and hence
\begin{eqnarray*}
\E(N_1 N_2) &  \leq & c_0^ 2 + o\left( n! \sum_{0 < k < 12
\ln n}
    (k+1) (n-k)! K_5 R(n,k) \right) \\
& & +\,n! \sum_{n - 2 \geq k > n- n^ {3/4} /2} (n-k+1) (2n\sp
{7/8})^ {n-k}
    K_5 R(n,k) \\
& & +\,n! \sum_{12 \ln n \leq k \leq n - 5e(n+3)^ {2/3}} 2 n^ 6
(n-k)!
    K_5 R(n,k) \\
& & +\,n! \sum_{12 \ln n \leq k \leq n - 5e(n+3)^ {2/3}} \lc 5 e
(n+3)^ {2/3} \rc
(2n^ {7/8})^ {\lc 5 e (n+3)^ {2/3} \rc - 1}
    K_5 R(n,k).
\end{eqnarray*}
It will be shown that the last three terms of the bound on $\E N_1^ 2$
are
respectively $O(1)$, $o(1)$, and $o(1)$.  This will show that $\E N_1^ 2
=
O(1)$, and also that the first four terms of the bound on $\E(N_1 N_2)$
sum to
$c_0^ 2 + o(1) +o(1) + o(1) = c_0^ 2 + o(1)$.  For large enough
$n$, the fifth
term of the bound on $\E(N_1 N_2)$ is bounded by
\begin{eqnarray*}
K_5 n!
&\times&n \lc 5 e (n+3)^ {2/3} \rc (2n^ {7/8})^ {\lc 5 e (n+3)\sp
{2/3} \rc - 1}
         R(n, n - \lc 5 e  (n+3)^ {2/3} \rc) \\
&   =  &(1 + o(1)) K_5 \left( \frac{5 e \pi}{2} \right)^ {1/2} n\sp
{11/24} \\
&      & \times
         \left( 1 + \frac{\lc 5 e (n+3)^ {2/3} \rc}{n} \right)\sp
         {-\left( n + \lc 5 e (n+3)^ {2/3} \rc \right)} \left(
\frac{8e}{n^ {1/8}}
         \right)^ {\lc 5 e (n+3)^ {2/3} \rc},
\end{eqnarray*}
which vanishes at a rate faster than any power of $n$.  Thus $\E(N_1 N\sb
2) =
c_0^ 2 + o(1)$, completing the proof of the theorem
via~(\ref{twobounds}).

The three estimates for the bound on $\E N_1^ 2$ are now routine
calculations.
Plugging in the value of $R(n,k)$ and using Stirling's
formula gives for the second term in the bound
\begin{eqnarray*}
& & (1 + o(1)) n! \sum_{0 < k < 12 \ln n} (k+1) (n-k)! K_5 R(n,k)
\\[2ex]
& = & (1+o(1)) K_5 n^ n e^ {-n} \sqrt{2 \pi n} \sum_{1 \leq k <
12 \ln n}
(k+1) (n-k)^ {n-k} e^ {-(n-k)}
\sqrt{2\pi(n-k)} \,  2^ {2n-2k} e^ {2n-k} \\
&&   \hspace{0.9in} \times (2n-k)^ {-(2n-k)} (n-k)^ {-1/2}
(2n-k)^ {-1/2} \\[2ex]
& = &(1 + o(1)) K_5 2^ {1/2} \pi \sum_{1 \leq k < 12 \ln n}
   (k+1) \left [ (n-k)^ {n-k} n^ n 2^ {2n-2k} (2n-k)^ {-(2n-k)}
\right ].
\end{eqnarray*}
Note that if the sum here were to contain a
$k=0$ term, that term would equal $1$.  Furthermore, changing $k$ to $k+1$
multiplies the part of the summand inside square brackets by
$$
(n-k-1)^ {-1} [(n-k-1)/(n-k)]^ {n-k} 2^ {-2}
[(2n-k-1)/(2n-k)]^ {-(2n-k)} (2n-k-1).
$$
Now $[(n-k-1)/(n-k)]^ {n-k} \to e^ {-1}$,
while $[(2n-k-1)/(2n-k)]^ {-(2n-k)} \to e$
and $(2n-k-1) / (n-k-1) \to 2$, all uniformly over
$k < 12 \ln n$.  Thus the successive ratios are $(1+o(1)) \frac{1}{2}$
uniformly over $k$ in the range of summation.  Therefore the sum is at most
$(1 + o(1)) \sum_{k=1}^ \infty (k+1) 2^ {-k} = (1 +o(1)) 3 = O(1)$,
establishing the first bound.

For the third term in the bound on $\E N_1^ 2$, let $m = n-k$.  Then
plugging in
for $R(n,k)$ and using Stirling's formula yields
\begin{eqnarray*}
&& n! \sum_{n - 2 \geq k > n- n^ {3/4} /2} (n-k+1) (2n^ {7/8})\sp
{n-k}
    K_5 R(n,k) \\[2ex]
& = &(1 + o(1)) K_5 \sum_{n-n^ {3/4}/2 < k \leq n-2} (n-k+1) (2n\sp
{7/8})^ {n-k}
   n^ n e^ {-n} \sqrt{2\pi n}\, 2^ {2n-2k} e^ {2n-k} \\
  &&\hspace{0.9in} \times (2n-k)^ {-(2n-k)} (n-k)^ {-1/2}
(2n-k)^ {-1/2} \\[2ex]
& \leq & (1 + o(1)) K_5 (2 \pi)^ {1/2} \sum_{2 \leq m < n\sp
{3/4}/2}
   (m^ {1/2}+ m^ {-1/2}) ((2e)^ {8/7}n)^ {7m/8} n^ n 2\sp
{2m}
 (n+m)^ {-(n+m)} \\[2ex]
& \leq & (1 + o(1)) K_5 (2 \pi)^ {1/2} \sum_{2 \leq m < n\sp
{3/4}/2}
   2m (8e)^ m n^ {-m/8} ,
\end{eqnarray*}
For the first inequality here we used
$\sqrt{n} (2n-k)^ {-1/2} = 1 + o(1)$ uniformly over $k > n - n\sp
{3/4}/2$,
and for the second we used $m^ {1/2} + m^ {-1/2} \leq 2m$ and $(n+m)\sp
{-(n+m)} \leq
n^ {-(n+m)}$.
Changing $m$ to $m+1$ multiplies the term by
$8e(1+1/m) n^ {-1/8}$, which vanishes in the limit
uniformly in $m$; thus the sum is dominated by the $m=2$ term,
whose value is a constant times $n^ {-1/4}$, and is thus $O(n\sp
{-1/4}) = o(1)$.

Finally, to bound the fourth term in the bound on $\E N_1^ 2$, plug in
to get
\begin{eqnarray*}
& &  n! \sum_{12 \ln n \leq k \leq n - 5 e (n+3)^ {2/3}} n^ 6
(n-k)!
      K_5 R(n,k)  \\[2ex]
&=& (1 + o(1)) 2 \pi K_5 \sum_{12 \ln n \leq k \leq n - 5e(n+3)\sp
{2/3}} n^ 6
     (n-k)^ {n-k} e^ {-(n-k)} \sqrt{n-k}\,n^ n e^ {-n} \sqrt{n}
2^ {2n-2k} e^ {2n-k}
     \\
&&   \hspace{1in} \times (2n-k)^ {-(2n-k)} (n-k)^ {-1/2} (2n-k)\sp
{-1/2} \\[2ex]
& \leq & (1 + o(1)) 2 \pi K_5 \sum_{12 \ln n \leq k \leq n - 5 e
(n+3)^ {2/3}}
   n^ 6 (n-k)^ {n-k} n^ n 2^ {2n-2k} (2n-k)^ {-(2n-k)},
\end{eqnarray*}
since $\sqrt{n} (2n-k)^ {-1/2} \leq 1$.

The sum here is at most $n$ times its largest term.
Let $r = k/n$ and rewrite the typical summand
as $n^ 6 [(4-4r)^ {1-r} / (2-r)^ {2-r}]^ n =
n^ 6 h(r)^ n$, say.  Now we find the maximum of $h(r)$ on $[0,1]$.
Taking logs gives
$$\ln h(r) = (1-r) \ln(4-4r) - (2-r) \ln(2-r),$$
so that
$$(d/dr) \ln h(r) = \ln(2-r) - \ln(4-4r) .$$
This increases from $-\ln 2$ to $\infty$ as $r$ increases
from $0$ to $1$, so the maximum of $\ln h(r)$ over the interval $(12
\ln n)/ n
\leq r \leq 1 - 5e(n+3)^ {2/3}/n$ is achieved at one of the endpoints, at
least
for large $n$.  Again for large enough $n$, we can, for any $\dd > 0$, get the
derivative of $\ln h$ on
$[0,(12 \ln n)/ n]$ to be bounded above by $-\ln 2 + \dd < -0.693 + \dd$, so
choosing $\dd = 0.003$ makes the derivative of $\ln h$ bounded above by
$-0.69$ on this interval.  Similarly, for large enough $n$ the derivative
on $[1-5e(n+3)^ {2/3}/n , 1]$ is bounded below by $1$.  Noting that
$h(0) = h(1) = 1$, it follows that the value of $n^ 6 h(r)^ n$ at
$r=(12 \ln n)/ n$ is at most $[e^ {(-0.69) (12 \ln n)/ n}]^ n n^ 6
= n^ {-2.28}$, and
the value at $1-5e(n+3)^ {2/3}/n$ is at most $[e^ {-5e(n+3)\sp
{2/3}/n}]^ n n^ 6 =
e^ {-5e(n+3)^ {2/3}} n^ 6 < n^ {-2.28}$ for large $n$.  Thus the
sum under
consideration is at most $n^ {-1.28}$.

Putting all of this together gives $\E N_1^ 2 \leq c_0^ 2 + c\sb
0 + O(1)
+ O(n^ {-1/4}) + O(n^ {-1.28}) = O(1)$, as desired.    $\Cox$

\section{Unoriented percolation}

In Section 5 we shall consider the first-passage time to $\1$ for unoriented
first-passage percolation on $\B_n$.  For completeness, in this section we
treat ordinary unoriented percolation and argue that the critical probability
is $1/n$, as put forth in the following theorem:

\begin{th} \label{unoriented percolation}
Let each edge of $\B_n$ be independently open with probability $p = c/n$,
$0 <
c < \infty$.  Then $\P(\0$ {\rm is connected to $\1$ by an (unoriented) open
path)}$\to (1 - x(c))^ 2$, where $x(c)$ is, as in Theorem~\ref{ordinary},
the
extinction probability for a Poisson($c$) Galton--Watson process.
\end{th}

\noindent{\bf Proof:}  Write $\theta_n \equiv \theta_n(c)$ for the
percolation probability in question.  We first note that $\limsup_n
\theta_n
\leq (1 - x(c))^ 2$ by a branching process approximation similar to that in
the second paragraph of the proof of Theorem~\ref{ordinary}; we omit the
details.

For the lower bound we may restrict attention to the case $c > 1$; it is
precisely for these values of $c$ that $y(c) := 1 - x(c) > 0$.  Let $0 < \ee
< y(c)/4$.  We rely heavily on a result of Ajtai, Koml\'{o}s, and
Szemer\'{e}di (1982):  $\P(F) \geq 1 - o(1)$, where $F$ is the event that
(a) there is exactly one component in the random graph formed by the open
edges
that has at least $(y(c) - \epsilon) 2^ n$ vertices, and (b) all the other
components are of size at most $\ee 2^ n$.

For $v \in B_n$, let $A_v$ denote the event $\{y$ is connected by an open
path to at least $(y(c) - 2\ee) 2^ n$ vertices$\}$.  By a simple application
of the FKG inequality (Fortuin, Ginibre, and Kasteleyn (1971)), the
indicators of the events $A_v$ are pairwise positively correlated.
Furthermore, conditionally given $F$, we have by symmetry
$$
\P(A_v) = \P(v \in \mbox{ the unique giant component (GC)}) = 2^ {-n}
\E(\mbox{size of GC}) \geq y(c) - \ee;
$$
thus, unconditionally, $\P(A_v) \geq y(c) - 2\ee$ for sufficiently large
$n$.  By FKG,
$$
\P(A_{\0} \cap A_{\1}) \geq (y(c) - 2\ee)^ 2,
$$
and so
\begin{eqnarray*}
\lefteqn{\P(\mbox{$\0$ and $\1$ are in the same component}\|F) \geq
  \P(\mbox{$\0$ and $\1$ are in the GC}\|F)} \\
 & & = \P(A_{\0} \cap A_{\1}\|F) \geq \frac{\P(A_{\0}
\cap A_{\1}) - (1 -
         \P(F))}{\P(F)} \geq (y(c) - 3\ee)^ 2,
\end{eqnarray*}
and hence $\theta_n \geq (y(c) - 4\ee)^ 2$, for sufficiently
large $n$.  Let
$\epsilon \downarrow 0$ to complete the proof.   $\Cox$

We close this section by noting that for $c < 1$ there is a more elementary
proof that $\theta_n(c) \to 0$.  For $x \in B_n$, let $g(x)$ denote the
probability that $x$ is connected to $\0$ by an open path.  Clearly, for $x
\not= \0$, $g(x)$ equals the probability that there is a neighbor $y$ of
$x$ such that $y$ is connected to $\0$ by an open path not containing $x$
and the edge $\{y,x\}$ is open.  Hence
\begin{equation} \label{gineq}
\mbox{$g(\0) = 1$, $g(x) \leq p \sum_{y \sim x} g(y)$ for $x \not=
\0$, $g(x)
\leq 1$ for all $x$,}
\end{equation}
where the sum is over vertices $y$ adjacent to $x$.

Repeatedly applying~(\ref{gineq}), we find $g(x) \leq pn = c$ for $x \not=
\0$, $g(x) \leq c^ 2$ for $d(\0,x) \geq 2$, $g(x) \leq c^ 3$ for
$d(\0,x)
\geq 3$, \ldots, and finally $\theta_n = g(\1) \leq c^ n \to 0$, as
desired.

\section{Richardson's growth model and unoriented percolation}

Consider the following model for the spread of disease.
Individuals are located at vertices of an
$n$-cube, with edges modelling pairs of individuals
in frequent contact.  One individual, $\0$, is infected
at time $0$ and the rest are healthy.  Independently for each edge
between an infected individual and an uninfected one,
there is a constant small probability per small unit of time
that the contact between those two individuals will cause
the uninfected one to become infected.
It is easy to construct from this description
a stochastic model for the growing set of infected individuals.
The model is a continuous time
Markov chain on the space of subsets of $\B_n$ which jumps
from $A$ to $A \cup \{ v \}$ at rate $k(A, v)$, where $k(A, v)$
is the number of infected neighbors of $v$, i.e., the number of neighbors of
$v$ in $A$.  This Markov chain is called Richardson's growth model.

Interesting questions about this model are (1) When
should we expect $\1$ to become infected? and (2) What
is the cover time, i.e., when should we expect
all the vertices to become infected?  In
this section we
discuss the first question, giving limiting upper and lower
of $1$ and $0.88$, respectively.  These are obtained by proving and then
exploiting the fact that the infection time for $\1$ in Richardson's model
has the same distribution as the first-passage time to $\1$ in {\em
unoriented\/} first-passage percolation on $\B_n$.
The cover time question is addressed in Section~6.

The following lemma reduces the problem of when $\1$ first becomes infected
to unoriented first-passage percolation with exponentially
distributed edge-passage times.
Since the oriented percolation
time is always at least as great as the unoriented percolation
time (the minimum over paths directed away from $\0$ must be at least the
minimum over all paths), the upper bound of Theorem~\ref{crossworld} for the
infection time of $\1$ is immediate.

\begin{lem} \label{reduction}
Let the edges $\{v,w\}$ of the undirected graph $\B_n$ be assigned
independent exponential random variables $X_{v,w} = X_{w,v}$ of mean
$1$.
Define the {\em infection time\/} of a vertex $v$, denoted $T_n(v)$, to
be
$\inf \sum_i X_{v_i v_{i+1}}$, where the $\inf$ is over all paths
from $v_0 =
\0$ to $v$.  Let $A(t) = \{ v \in \B_n \,:\, T_n (v) \leq t \}$.  Then
the random map $A$ from $[0 , \infty )$ to subsets of $\B_n$
has the same law as Richardson's model.
\end{lem}

\noindent{Proof:}  See Durrett (1988, page 177) for a sketch of this proof.   
$\Cox$

\begin{th} \label{crossworld}
For any $\ee > 0$, the probability of finding $\1$ infected
by time $1+\ee$ in Richardson's model on $\B_n$, beginning with
only $\0$ infected at time $0$, tends to $1$ as $n \rightarrow \infty$.
\end{th}

\noindent{\bf Proof:}  Theorem~\ref{OFPP time} and Lemma~\ref{reduction}.
$\Cox$

We doubt whether this result is sharp, since there is no
reason why the unoriented percolation time should be
as great as the oriented percolation time.  The next theorem,
based on a calculation by R.~Durrett (personal communication), gets
a lower bound for the unoriented first-passage time by comparing
to a branching translation process (BTP).  This is a process,
started with a single particle at $\0$, for which each existing particle
generates offspring at rate $n$, where the offspring are each
displaced from the parent by an independent uniform random step $e(j)$.
Letting $Z(x,t)$ be the number of particles at $x$ at time $t$,
the process $(Z(x,t) : x \in \B_n, \, t \geq 0)$ is formally
defined by the transition rates $(Z(x)) \rightarrow (Z(x)
+ \dd_{xy})$ at rate $\sum_{w:\,d(y,w) = 1} Z(w)$, where $\dd_{xy}$
is $1$ if $x=y$ and $0$ otherwise.  It is easy to couple BTP to
Richardson's model so that the set of infected vertices in Richardson's
model is always
a subset of the set of populated vertices in BTP.  Thus the
first time $\tau_n$ that $\1$ is populated in BTP is stochastically
less than the first infection time $T_n$ of $\1$ in Richardson's model.

\begin{th}[Durrett] \label{lowerDurrett}
As $n \rightarrow \infty$, the time $\tau_n$ of first population
of $\1$ in BTP converges in probability to $\ln(1+\sqrt{2})
\doteq 0.88$.  Consequently, $\P (T_n \leq \ln(1+\sqrt{2}) - \ee)
\rightarrow 0$.
\end{th}

\noindent{\bf Proof:}  The lower bound will be gotten by a
routine first moment calculation.  The upper bound in probability
for BTP (which is not necessary for the result on Richardson's model) requires
a second moment calculation and a little more work.  Fix $n$ and write
$m_1(x, t)$ for $\E Z(x,t)=$ the expected number of particles
at $x$ at time $t$ in BTP starting from a single particle at $\0$.
We remark for later that this is also the expected number of
offspring at $y+x$ at time $s+t$ of a particle at $y$ at time $s$
that are born to the particle after time $s$, where the addition in $y+x$
is taken, as usual, to be coordinatewise mod~$2$ addition.
Since $\P (Z(\1 , t) > 0) \leq m_1(\1 , t)$, the lower bound in probability
will follow from showing that $m_1(\1 , t) \rightarrow 0$ as $n \rightarrow
\infty$ for any $t < \ln(1+\sqrt{2})$.  Viewing vertices of $\B_n$ as
sets, we
write $|x|$ for the cardinality of $x$; the differential equation for
$m_1(x, t)$ is easily seen to be
\begin{equation} \label{diffeq}
{d\,m(x,t) \over dt} = \sum_{y:\,d(x,y)=1} m(y,t)
\end{equation}
with initial conditions $m(x,0) = \dd_{\0,x}$.  Let
$$
p(x, t) = \left ( {1 - e^ {-2t} \over 2} \right )^ {|x|}
    \left ( {1 + e^ {-2t} \over 2} \right )^ {n - |x|}
$$
be the probability that a simple random walk with rate $n$ started
at $\0$ is at $x$ at time $t$.  Then, as may be verified by a variety of
probabilistic and analytic arguments, the unique solution to~(\ref{diffeq}) is
given by
$$m_1(x, t) = e^ {nt} p(x, t).
$$
Putting $x=\1$ gives
$$m_1(\1 , t) = \left ( {e^ t - e^ {-t} \over 2} \right )^ n .$$
Since $(e^ t - e^ {-t})/2$ is increasing in $t$ and equal to $1$ at
$t= \ln(1+\sqrt{2})$, it follows that for $t < \ln(1+\sqrt{2})$,
$m_1(\1 , t)$ tends to $0$ as $n \rightarrow \infty$.  Hence
$\P (\tau_n \leq t) \rightarrow 0$, as desired.

The upper bound in probability on $\tau_n$ is gotten by a now familiar
sort of
argument.  Fix $\ee > 0$.  First
the second moment method
is used to show that $\liminf_n \P (\tau_n \leq \ln(1+\sqrt{2}) +
\ee) \geq
1/12$.  Then the initial
branching of the process is used to show that with just $2\ee$ more time units,
there are actually many independent chances of no less than $1/12$ each to get
$\1$ populated, and hence the probability that this occurs is near $1$.

Begin with
$$
\P (\tau_n \leq t ) = \P (Z (\1,t) > 0) \geq {(\E\,Z(\1 ,t))^ 2
    \over \E (Z (\1 , t)^ 2)} = {m_1(\1 , t)^ 2 \over m_2(\1, t)},
$$
where $m_2(x, t) := \E (Z(x,t)^ 2)$.  To calculate the value of $m\sb
2(x, t)$
in terms of $m_1(x, t)$, write $Z(x,t)^ 2$ as $Z(x,t)$ plus twice the
number of unordered pairs of distinct particles at $x$ at time $t$.
Each such pair of particles has a well-defined time $s$ at which their
ancestral lines first split apart.  Say that at time $s$ a particle $p_1$
at
vertex $y$ gave birth to a particle $p_2$ at vertex $y + e(i)$, and that
both
particles are descendants of $p_1$ but only one is a descendant of $p\sb
2$.
For fixed $y$ and $i$ and interval $[s, s+ds)$, the expected number of
such pairs is $m_1(y, s)\,ds\,m_1(x - y, t - s)\,m_1(x - y - e(i),
t - s)$,
so summing over $y$ and $i$ and integrating over $s$ gives
\begin{equation} \label{msub2}
{m_2(\1 , t) \over m_1(\1, t)^ 2} = {1 \over m_1(\1 , t)} +
\sum_i 2
    \int_0^ t\,ds \sum_y {m_1(y, s) m_1(\1 -y, t-s)
m_1(\1 -
y - e(i), t - s)
    \over m_1(\1 , t)^ 2} .
\end{equation}
Now fix $t = \ln(1+\sqrt{2}) + \ee$.  The first term tends to $0$ as
$n \rightarrow \infty$, so it suffices to show that the $\limsup$
of the sum on $i$ is at most $12$.

Substituting $m_1(x, t) = e^ {nt} p(x, t)$ into the sum on $i$
in~(\ref{msub2})
yields
\begin{equation} \label{integrand}
2 \int_0^ t\,ds \sum_{y,i} e^ {-ns} {p(y, s) p(\1 - y, t - s)
p(\1 - y -
 e(i), t - s) \over p(\1 , t)^ 2}.
\end{equation}
Next, plug in the value for $p(x, t)$.  At the same time, group together
all $y$ on the same level of $\B_n$, i.e., all $y$
with $|y| = k$ for each $k$.  Then $|y+e(i)|$
will equal either $k+1$ or $k-1$; since $p(\1 - x, t)$ increases with $|x|$,
we get an upper bound by replacing $|y+e(i)|$ by $k+1$.  This gives an upper
bound for the integrand of
\begin{eqnarray*}
&& n e^ {-ns} {1 + e^ {-2(t-s)} \over 1 - e^ {-2 (t-s)}} \\
   &&\times \sum_{k=0}^ n 2^ {-n} {n \choose k} [(1-e\sp
{-2s})(1+e^ {-2(t-s)})^ 2]^ k
      [(1+e^ {-2s})(1-e^ {-2(t-s)})^ 2]^ {n-k}
(1-e^ {-2t})^ {-2n}.
\end{eqnarray*}
The sum over $k$ is just the binomial expansion of
$$
\left ( { (1-e^ {-2s})(1+e^ {-2(t-s)})^ 2 +
 (1+e^ {-2s})(1-e^ {-2(t-s)})^ 2
   \over 2(1-e^ {-2t})^ 2} \right )^ n,
$$
and simplifying this yields
$$
\left( {1 \over 2(1-e^ {-2t})^ 2} (2 + 2e^ {-4(t-s)} -4
e^ {-2t}) \right)^ n
 = \left( 1 + {e^ {-4 (t-s)} - e^ {-4t} \over (1 - e\sp
{-2t})^ 2} \right)^ n,
$$
which gives a bound for the integrand in~(\ref{integrand}) of
\begin{equation} \label{sloppy bound}
n \; {1 + e^ {-2(t-s)} \over 1 - e^ {-2 (t-s)}} e^ {-ns}
\left ( 1 + {e^ {-4 (t-s)} - e^ {-4t} \over (1 - e^ {-2t})^ 2 }
\right )^ n.
\end{equation}

We need a better bound on the integrand when
$s$ is near $t$: the factor $\displaystyle{{1 + e^ {-2(t-s)} \over 1 -
e^ {-2
(t-s)}}}$ blows up like $(t-s)^ {-1}$, which is not integrable.  Note
that in
the case $k=n$ it is not possible to have $|y+e(i)| = k+1$. Thus
for the $k=n$ term, the factor $\displaystyle{{1 + e^ {-2(t-s)} \over
1 - e^ {-2
(t-s)}}}$ can be replaced by its reciprocal.  This reduces the integrand
significantly when the $k=n$ term is the dominant term in the sum.
The ratio of the $k=n$ term to the entire above sum on $k$ is
\begin{eqnarray*}
&& [(1-e^ {-2s})(1+e^ {-2(t-s)})^ 2]^ n /
 [(1-e^ {-2s})(1+e^ {-2(t-s)})^ 2 +
   (1+e^ {-2s})(1-e^ {-2(t-s)})^ 2]^ n \\[2ex]
& \geq & [(1-e^ {-2s}) / ((1-e^ {-2s})
+ (1+e^ {-2s})(2(t-s))^ 2)]^ n
 \\[2ex]
& \geq & [1 - K (t-s)^ 2]^ n \\[2ex]
& \geq & 1 - K(t-s)
\end{eqnarray*}
for some constant $K$ when $t - s \leq 1/n$.
($K = 15$ will do when $n \geq 2$.)
Also, $\displaystyle{[{1 + e^ {-2(t-s)} \over
1 - e^ {-2 (t-s)}}]^ {-2} = (t-s)^ 2 + O((t-s)^ 3)}$.
Putting this all together,
a better bound for the integrand in~(\ref{integrand}),
uniformly for $s$ satisfying $t - s \leq 1/n$, is
\begin{equation} \label{better bound}
(1 + o(1)) K (t-s) n {1 + e^ {-2(t-s)} \over 1 - e^ {-2 (t-s)}} e\sp
{-ns}
\left ( 1 + {e^ {-4 (t-s)} - e^ {-4t} \over (1 - e^ {-2t})^ 2 }
\right )^ n.
\end{equation}

In order to make use of~(\ref{sloppy bound}) and~(\ref{better bound}),
examine the function
$$G(s,u) = \ln \left[ e^ {-s}  \left( 1 + {e^ {-4 (u-s)} -
e^ {-4u} \over (1 - e^ {-2u})^ 2 } \right) \right] .$$
This is convex in $s$ for the values of $u$ we are interested in, which
may be seen by differentiating twice with respect to $s$:
writing $c \equiv
c_u
= e^ {-4u} / (1-e^ {-2u})^ 2$ gives
\begin{eqnarray*}
{\partial^ 2 G(s,u) \over \partial s^ 2} & = & {\partial^ 2
\over
\partial s^ 2} \left [ -s + \ln(1 + c(e^ {4s} - 1)) \right ] \\[2ex]
& = & {\partial \over \partial s} \left ( -1 + {4ce^ {4s} \over
    1+ c(e^ {4s} - 1)} \right ) \\[2ex]
& = & \frac{16 c (1 - c) e^ {4s}}{[1 + c(e^ {4s} - 1)]^ 2},
\end{eqnarray*}
which is positive for all $s$ whenever $c < 1$.  This is true if and
only if $u > \ln\sqrt{2}$
and hence for $u = \ln(1+\sqrt{2})+\ee$.
In particular, the maximum of $G(s,u)$ over $s \in [0,u]$ is achieved at an
endpoint. But $G(0,u) = 0$ and
$$
G(u,u) = \ln \left[ e^ {-u} (1 + {1-e^ {-4u} \over (1-e^ {-2u})\sp
2}) \right]
       = \ln \left[2 e^ u / (e^ {2u} - 1) \right].
$$
This is decreasing in $u$ and has value $0$ when $u = \ln(1+\sqrt{2})$.  Thus
when $u = t = \ln(1+\sqrt{2}) + \ee$, $G(u,u)$ has a negative value which we
shall call $-V(\ee)$.  Now we bound the upper bound for $(m_2(\1 , t) -
m_1(\1, t)) / m_1(\1, t)^ 2$ given by display~(\ref{integrand}) in
three pieces:
\begin{eqnarray}
&& [m_2(\1 , t) - m_1(\1 , t)] / m_1(\1, t)^ 2 \nonumber
\\[2ex]
& \leq & 2\int_0^ {1/2}\,{1 + e^ {-2(t-s)} \over 1 - e\sp
{-2 (t-s)}}
   n e^ {n G(s,t)}\,ds \label{first} \\[2ex]
& + & 2\int_{1/2}^ {t-1/n}\,{1 + e^ {-2(t-s)} \over 1 - e\sp
{-2 (t-s)}}
   n e^ {n G(s,t)}\,ds \label{second}\\[2ex]
& + & (1 + o(1)) 2\int_{t-1/n}^ t\,{1 + e^ {-2(t-s)} \over 1 -
e^ {-2 (t-s)}}
   n K(t-s) e^ {n G(s,t)} \label{third}\,ds .
\end{eqnarray}

For the first piece we calculate the
value of $G(1/2, t)$, getting a constant less than $-1/4$.
Thus by convexity, $G(s, t) \leq -s/2$ for $0 \leq s \leq 1/2$.
Now for $0 \leq s \leq 1/2$ and any $\ee$, the factor $\displaystyle{{1 +
 e^ {-2(t-s)}
\over 1 - e^ {-2 (t-s)}}}$ is at most $3$, so the contribution
from~(\ref{first})
is at most $2\int_0^ {1/2}\,3n e^ {-ns/2}\,ds = 12(1 -
e^ {-n/4}) < 12$.

For the second piece, bound $G(s,t)$ on $1/2 \leq s \leq t$
by its values at the endpoints; for small $\ee$ the greater value is
the value at the right endpoint, namely, $-V(\ee)$.  The value
of $\displaystyle{{1 + e^ {-2(t-s)} \over 1 - e^ {-2 (t-s)}}}$ on
$1/2 \leq s
 \leq
t-1/n$ is bounded by its maximum, which is achieved at $s=t-1/n$ and
has a value of at most $n + 1$.  Thus the contribution from~(\ref{second})
is at most $2 \int_{1/2}^ {t-1/n}\,(n+1) n e^ {-nV(\ee)}\,ds < 2 t
(n + 1) n
e^ {-n V(\ee)}$, and this tends to $0$ as $n \rightarrow \infty$.

To bound the third piece, expand the integrand in powers of $(t-s)$
to compute the integral as $(1+o(1)) \int_{t-1/n}^ t\,(t-s)^ {-1}
n K(t-s) e^ {n
G(s,t)}\,ds \leq (1+o(1)) K e^ {-nV(\ee)}$, which tends to $0$ as
$n \rightarrow \infty$.  Thus the entire integral is bounded by $12
(1 + o(1))$,
and $\liminf_n \P (\tau_n \leq \ln(1+ \sqrt{2}) + \ee) \geq 1/12$.

Finally, to show that $\P (\tau_n \leq \ln(1+\sqrt{2}) + 3 \ee)
\rightarrow 1$, let $A$ be the set of particles at distance $1$ from
$\0$ at time $\ee$.  Then $|A| > n^ {1/2}$ with probability
approaching $1$ as $n \rightarrow \infty$.  The offspring of
elements of $A$ now act independently from time $\ee$ to time
$\ln(1+\sqrt{2}) + 2\ee$, each particle at some $y$ with
$d(\0,y)=1$ having probability at least $(1 + o(1)) \frac{1}{12}$ of having a
descendant at the antipodal point to $y$ at time
$\ln(1+\sqrt{2}) + 2\ee$, according to the calculation
just completed.  Letting $B$ be the set of particles
at sites that neighbor $\1$ at time $\ln(1+\sqrt{2}) + 2\ee$,
$\P (|B| > n^ {1/4})$ is at least $\P(|A| > n^ {1/2})
\P (X > n^ {1/4})$, where $X$ is a binomial with parameters
$\lc n^ {1/2} \rc$ and $(1 + o(1)) \frac{1}{12}$.  This probability also
tends
to $1$ as $n \rightarrow \infty$.  Finally, $\P (\tau_n \leq
\ln(1+\sqrt{2}) +
3\ee \||B| > n^ {1/4}) \geq 1 - e^ {-\ee n^ {1/4}}$,
which tends to $1$ as
$n \rightarrow \infty$, proving the theorem.    $\Cox$

\section{Covering times in Richardson's model}

This section answers affirmatively the question of whether
the time until the entire $n$-cube $\B_n$ is infected is
bounded in probability as $n \rightarrow \infty$.  The
constant upper bound given here is $4 \ln(4 + 2 \sqrt{3}) + 6 \doteq 14.04$.
This is by no means sharp, but on the other hand we also produce a lower
bound
on the cover time of $\frac{1}{2} \ln(2+\sqrt{5}) + \ln 2 \doteq 1.41$.
Since the infection time of $\1$ is bounded between $\ln(1+\sqrt{2}) \doteq
0.88$ and $1$ in probability, this means that the $\liminf$ in probability
of the cover time is strictly greater than the $\limsup$ in probability of
the time to reach the farthest vertex.  Perhaps the cover time has a
limit in probability, but we do not venture a guess as to
what the limit should be.

\subsection{The upper bound}

The following statement of duality in Richardson's model will be
helpful.  The proof can be found in any introduction to the contact
process, such as Durrett (1988).  The intuition is to think of
$A_2 (t-s)$ as the set of vertices that would be able to
infect $\1$ by time $t$ if they were infected at time $s$.

\begin{lem} \label{duality}
Let $(A_1(t))$ be Richardson's model on $\B_n$ as defined
above and let $(A_2(t))$ be an independent copy with the
difference that $A_2(0)$ is set to be $\{ \1 \}$ instead
of $\{ \0 \}$.  Then $\P(\1 \in A_1(t)) = \P(A_1(s)
\cap A_2(t-s) \neq \emptyset)$ for any $0 \leq s \leq t$.   $\Cox$
\end{lem}

\begin{th} \label{highprob}
Let $(A(t))$ be Richardson's model on $\B_n$ starting with
only $\0$ infected at time $0$.  Let $c := 2 \ln(4 + 2 \sqrt{3}) + 3 \doteq
7.02$.  Then for any $\ee > 0$ there
exists $N = N(\ee)$ and $h(\ee) > 0$ such that for $n \geq N$,
$\P(\1 \in A(c + \ee)) \geq 1 - (4 + h(\ee))^ {-n}$.
\end{th}

We defer the proof of Theorem~\ref{highprob} in order to present the resulting
cover time upper bound.

\begin{cor} \label{cover}
For any $\ee > 0$, $\P (A (2c+\ee ) = \B_n ) \rightarrow 1$
as $n \rightarrow \infty$.
\end{cor}

\noindent{{\bf Proof} of} Corollary~\ref{cover}:  Let $N \equiv N(\ee / 2)$ be
chosen as in the statement of
the previous theorem and pick any $n \geq 2N$.  Let
$x$ be any element of $\B_n$.  First suppose that the distance
from $\0$ to $x$ is $d \geq n/2$.  Then the sublattice
with top element $x$ and bottom element $\0$ is a Boolean algebra of rank at
least $N$, and the induced process on the sublattice (i.e., the process with no
infections allowed except on the sublattice) is still Richardson's model.
By the previous theorem, $x$ is infected
by time $c+ \ee / 2$ with probability at least
$1 - (4 + h(\ee/2))^ {-d} \geq 1 - (4 + h(\ee/2))^ {-n/2}$.

On the other hand, suppose $x$ is at distance less than $n/2$ to $\0$.  Then
after time $c + \ee / 2$ the top element $\1$ is infected with
probability at least $1-(4 + h(\ee/2))^ {-n}$ and, conditioned on that, the
probability that $x$ is infected another $c + \ee / 2$
time units later is at least $1 - (4 + h(\ee/2))^ {-n/2}$ by the previous
argument (since the distance from $x$ to $\1$ is at least $n/2$).
Thus each $x$ fails to be infected at time $2c + \ee$ with
probability at most $(4 + h(\ee/2))^ {-n} + (4 + h(\ee/2))^ {-n/2}$, and
summing
over all $x$ gives at most $(1 + h(\ee/2)/4)^ {-n/2} + 2^ {-n}$, which
tends
to $0$ as $n \rightarrow \infty$.   $\Cox$

Notice that the reason we get $2c$ instead of $c$
as an upper bound in probability for the cover time is that
the proof of Corollary~\ref{cover} gives better upper bounds on the probability
that a vertex is uninfected the further it is from $\1$.
It is unlikely that the bounds reflect the true state of affairs.
In particular, we suspect that the random time until a vertex
$x$ is infected is stochastically increasing in $|x|$.  This would
immediately imply an upper bound in probability of $c \doteq 7.02$ for the
covering time.
(In fact, the bound could then be lowered to $\tilde{c} := 2 \ln(2 +
\sqrt{2}) + 3 \doteq 5.46$ by establishing the modification $\P(\1 \in
A(\tilde{c} + \ee)) \geq 1 - (2 + \tilde{h}(\ee))^ {-n}$ of
Theorem~\ref{highprob}.)
More generally, we have the following conjecture.
\begin{conj} \label{stoch}
Let $G$ be a graph with distinguished vertex $x$.  Let $H$ be the
graph $G \times \{ 0,1 \}$, with edges between $(x,i)$ and $(y,i)$
for neighbors $x,y$ of $G$ and $i=0,1$ and edges between $(y,0)$ and
$(y,1)$ for all $y \in G$.  For $z \in H$, let $T(z)$ be the
time that $z$ is first infected in Richardson's model on $H$
beginning with a single infection at $(x,0)$.  Then, for any $y \in G$,
$T(y,0)$ is stochastically smaller than $T(y,1)$.
\end{conj}

We now turn to the proof of Theorem~\ref{highprob}.

\noindent{\bf Proof} of Theorem~\ref{highprob}:
Throughout the proof we use the notation
$x'$ for the complement of a vertex $x$ (viewing $x$ as a subset of $\{
1,\ldots, n \}$) and $S'$ for $\{x' \in \B_n: x \in S\}$ when $S \subseteq
\B_n$. Let $t_0 = 0$, $t_1 =
\ln(4 + 2 \sqrt{3}) + \ee / 7$, $t_2 = \ln(4 + 2 \sqrt{3}) + 2 \ee / 7$,
$t_3
= \ln(4 + 2 \sqrt{3}) + 1 + 3 \ee / 7$, and $t_4 = \ln(4 + 2 \sqrt{3}) +
2
+ 4 \ee / 7$.  By Lemma~\ref{duality} it suffices to find $h$ for which
$\P(A_1(t_4) \cap A_2(t_3) \neq \emptyset) \geq
1-(4+h(\ee))^ {-n}$ for large
$n$. The method will be to watch the evolutions of $A_1$ and $A_2$ and
look for
vertices $x$ for which simultaneously $x \in A_1(t)$ and $x' \in
A\sb2(t)$.  In particular, we will show that for $i=1,2,3$, certain subsets
$D_i$ of $\{ x : x \in A_1(t_i) \mbox{ and } x'
\in A_2(t_i) \}$
%*% What I've written on the next line is a new correction.
(actually, of a slight modification of this for $D_3$)
are sufficiently large, and
then we will argue that each $x$ in $D_3$ has an independent chance of
becoming
an element of $A_1(t_4) \cap A_2(t_3)$.

For $1 \leq i,j \leq n$ let $\eta_{ij}  = \{ i,j \} \in \B_n$.
Let $n_1 = \lfloor n - 2 \log_2 n \rfloor$ and let
$S_2 \subseteq \B_n$ be the set $\{ \eta_{ij} \, : \, i < j
\leq n_1 \}$.  Let $S_1 = \{ \{ i \} : i \leq n_1 \}$ be the
set of elements at level $1$ of $\B_n$ beneath $S_2$.  For each $x \in
S_2$ let
$T_x \subseteq \B_n$ be the set $\{\,y : x \cap \{ 1 , \ldots ,
n_1 \} = y \cap
\{ 1, \ldots ,  n_1 \}\,\}$.  Note that $T_x$ and
$T_y$ are disjoint for distinct $x,y \in S_2$; similarly for
$T_x'$ and $T_y'$.

Several of the arguments below will involve the monotonicity
of Richardson's model: forbidding some edges to pass
the infection at various times decreases $A(t)$ and hence can only increase
the
infection time to any vertex. It can therefore only increase all infection
times to suppose for $A_1$ that from time $t_0$ to time $t_1$
infections
may occur only in rank $1$ of $\B_n$, from time $t_1$ to time $t_2$
infections
may occur only in rank $2$ (and in fact later the set of allowed infections
will be further restricted), and from time $t_2$ to time $t_3$
infections may
occur only between $x$ and $y$ when the symmetric difference $x
\bigtriangleup
y$ is the singleton $\{ i \}$ for some $i > n_1$ (in other words,
infections
from $t_2$ to $t_3$ may occur only when $x$ and $y$ are neighbors both
in $T_z$
for some $z$). Also suppose dually for $A_2$ that infections
between times $t_0$ and $t_1$ occur only in level $n-1$,
that infections between times $t_1$ and $t_2$
occur only in level $n-2$, and that infections
between times $t_2$ and $t_3$
occur only between elements of the same $T_z'$.  Finally, suppose
that between times $t_3$ and $t_4$, infections in $A_1$ spread only
between
vertices $x$ and $y$ for which $x \cap \{ n_1 + 1 , \ldots , n\}
= y \cap \{ n_1 + 1 , \ldots , n\}$.

Let $U_1 = S_1 \cap A_1 (t_1)$, let $V_1 = S_1'
\cap A_2 (t_1)$,
and let $D_1 = U_1 \cap V_1'$.  The first claim is that
there exist $\dd_1 (\ee )$ and $h_1 (\ee )$, both positive,
for which $\P (|D_1| \leq \dd_1(\ee) n)
\leq (4+h_1(\ee) )^ {-n}$ for
sufficiently large $n$.  This is just a large deviation
calculation.  $|D_1|$ is the sum of $n_1$ i.i.d.\ Bernoulli
random variables, each equalling $1$ with probability
$p = (1-e^ {-t_1})^ 2$ and $0$ with probability $q =
2e^ {-t_1} - e^ {-2 t_1}$.
Let $0 < a < p$ and write $b = 1 - a$.  By choosing the optimal value $\ln
(\frac{bp}{aq})$ for $\theta > 0$, we find from the moment generating function
inequality
$$
\P(|D_1| \leq a n_1) = \P(e^ {-\theta |D_1|} \geq e^ {-\theta
a n_1}) \leq
e^ {\theta a n_1} \E e^ {-\theta |D_1|}
$$
that $\P(|D_1| \leq a n_1)
\leq [(p/a)^ a (q/b)^ b]^ {n_1}$.  As $a \downarrow
0$, $(p/a)^ a (q/b)^ b \downarrow q = 2e^ {-t_1} -
e^ {-2 t_1} <
2 (4+2\sqrt{3})^ {-1} - (4+2\sqrt{3})^ {-2} = (2 - \sqrt{3}) -
(7/4 - \sqrt{3})
= 1/4$.  Thus there exists $\dd_0(\ee) > 0$ and $h_0(\ee) > 0$ such
that $\P(|D_1| \leq \dd_0(\ee) n_1)
\leq (4+h_0(\ee))^ {-n_1}$.
Since $n_1 / n \rightarrow 1$, this implies the existence of
$h_1(\ee)$
such that for any fixed $\dd_1 < \dd_0(\ee)$,  $\P(|D_1| \leq
\dd_1 n) \leq
(4+h_1(\ee))^ {-n}$ for sufficiently large $n$.

Now let $U_2 = S_2 \cap A_1 (t_2)$, let $V_2 = S_2'
\cap A_2
 (t_2)$,
and let $D_2 = U_2 \cap V_2'$.  The next claim is that there
exist
$\dd_2 (\ee )$ and $h_2 (\ee )$, both positive,
for which $\P (|D_2| \leq \dd_2 n^ 2) \leq (4+h_2(\ee) )^ {-n}$
for
sufficiently large $n$.  Between times $t_1$ and $t_2$ the only
infections
allowed in $A_1$ involve vertices in $S_1$ infecting vertices in
$S_2$.  For
convenience, restrict further the allowed infections
by requiring that $\{ i \}$ may infect
$\{ i,j\}$ only if $i < j$ and $j - i < n_1 / 2$ or $i > j$ and
$i - j > n_1 /
2$.  Note that each $\{ i,j \}$ in $S_2$ can be infected by $i$ or $j$
but not
both. The exception is $|j - i| = n_1 / 2$; in that case $\{i, j\}$ cannot
be infected at all.  Then for each $x \in D_1$, there is a set $S(x)$ of
$\lc n_1 / 2 \rc - 1$ vertices in $S_2$
that $x$ can infect between times $t_1$ and $t_2$, and these
sets are disjoint as $x$ varies over $S_1$.  For each $x \in D_1$,
the number of vertices infected by $x$ by time $t_2$ whose complements
have been infected by $x'$ by time $t_2$ in the process $A_2$ (with
the dual
restrictions) is a binomial random variable with parameters $\lc n_1 / 2
\rc -
1$ and $(1 - e^ {-\ee / 7})^ 2$.  Now the probability that this
binomial is less
than half its mean is exponentially small in $n$,
say $\leq e^ {-\alpha n}$, so
conditioning on $|D_1| > \dd_1 n$, the probability that no more than
$\lc \dd_1
n \rc / 2$ of these i.i.d.\ binomials are greater than half their means is at
most $2^ {\lc \dd_1 n \rc} e^ {- \alpha \dd_1 n^ 2 / 2}$.
This is smaller than $(4+h_1 (\ee ))^ {-n}$ for sufficiently large $n$,
and
thus we have shown that $\P (|D_2| \leq (1 - e^ {-\ee / 7})^ 2 \lc
\dd_1 n \rc
\left( \lc n_1 / 2 \rc - 1 \right) / 4) < 2 (4+h_1 (\ee ))^ {-n}$
for
sufficiently large $n$.  Choosing $\dd_2 < (1 - e^ {-\ee / 7})^ 2
\dd_1 / 8$
and $h_2 < h_1$ proves the second claim.

Condition until the last sentence of this paragraph on $D_2$.
Between times $t_2$ and $t_3$ the spread of infection in $A_1$ is
confined
to each $T_z$, so the spread of infection is independent on each $T_z$.
The same goes for the propagation of $A_2$ on each $T_z'$.
On each $T_z$, the process is just a Richardson's model
on a cube of dimension $n-n_1$; hence uniformly for $x \in D_2$ and
$y \in T_x$
%*% I wrote you e-mail about the need for the next line.
with $\frac{1}{4} (n - n_1) \leq d(x,y) \leq \frac{3}{4} (n - n_1)$,
the probability that $y$ is infected by time $t_3 = t_2 + (1 + \ee/7)$
tends to $1$ for large $n$ by Theorem~\ref{crossworld}.
Similarly, the probability that $\hat{y}$ is dual infected by time
$t_3$ tends uniformly to $1$, where $\hat{y}$ is the element of
$T_x'$ that agrees with $y$ in the last $n-n_1$ places.
%*% I wrote you e-mail about the need for the next line.
(Note $d(x',\hat{y}) \geq \frac{1}{4} (n - n_1)$.)
Thus the probability that $y$ is infected by $t_3$ and
$\hat{y}$ is dual infected by $t_3$ also tends uniformly to $1$.
Thus uniformly for $x \in D_2$, the expected cardinality of $V_x :=
\{ y \in T_x : y \mbox{ is $A_1$-infected by }t_3
\mbox{ and } \hat{y} \mbox{ is $A_2$-infected by } t_3 \}$
is at least $(1-o(1)) |T_x| \geq
(1-o(1)) n^ 2$.  Then, since we always have $|V_x| \leq |T_x|$,
$\P (|V_x| \leq
n^ 2 / 2) \leq 1/2$ for large $n$,
so $\P (|V_x| \leq n^ 2 / 2 \mbox{ for all } x
\in D_2 ) \leq (1/2)^ {|D_2|}$.  Combining this with the previous
claim
about the distribution of $|D_2|$
shows that the (now unconditional) probability of the event $F$ that there is
some $x \in D_2$ for which $|V_x| \geq n^ 2 / 2$ is at least $1 -
(4+h_2 (\ee
))^ {-n} - (1/2)^ {\dd_2 n^ 2}
\geq 1 - (4+h_3 (\ee ))^ {-n}$ for $h_3 < h_2$ and
sufficiently large $n$.

Finally, condition on $F$.  Let $x$ be an element of $D_2$
with $|V_x| \geq n^ 2 / 2$ and let $D_3 = V_x$.
It suffices to show that
with high probability there is some $y \in D_3$ such that
$\hat{y}$ is infected at time $t_4$, since we already know
that $\hat{y}$ is dual infected at time $t_3$.  Recall that the
only way that infection spreads between times $t_3$ and $t_4$
is between neighboring vertices that intersect $\{ n_1 + 1 ,
\ldots, n \}$
equally.  This effectively breaks $\B_n$ into $2^ {n-n_1}$
fibers on which
propagation of the infection is independent and behaves like a Richardson's
model of dimension $n_1$.  Each element of $V_x$ is in a different
fiber,
and uniformly for $y \in V_x$ the probability that infection passes
from $y$ to the opposite corner $\hat{y}$ of the fiber in the time $t_4 -
t_3 =
1 + \ee/7$ tends to $1$ for large $n$ by
Theorem~\ref{crossworld} and in particular is eventually greater that
$1/2$.
Thus the conditional probability given $F$ that some $\hat{y}$ is infected at
time $t_4$ is at least $1- (1/2)^ {n^ 2 / 2}$ for large enough $n$.

Combining all the conditional probabilities and using
Lemma~\ref{duality} yields a probability
of at least $1 - (4+h_3 (\ee ))^ {-n} - (1/2)^ {n^ 2 / 2}$
that $\1$ is infected by time $t_3 + t_4 = 2\ln(4+2\sqrt{3}) + 3 +
\ee$.
Pick $h(\ee) < h_3 (\ee)$.  Then when $n$ is large enough,
our bound is at least $1 - (4+h (\ee ))^ {-n}$ and the proof is finsihed.
$\Cox$

\subsection{The lower bound}

\begin{th} \label{separation}
For any $\ee > 0$, $\P (A(\frac{1}{2} \ln(2 + \sqrt{5}) + \ln 2 - \ee) =
\B_n)
\rightarrow 0$ as $n \rightarrow \infty$.
\end{th}

\noindent{\bf Proof:}  First stochastically dominate Richardson's model at
time
$t = \frac{1}{2} \ln(2+\sqrt{5}) - \ee / 2$ by the corresponding value of a
branching translation process, as in Section~5.  Let $x$ be any vertex with
$d(\0, x) \geq n/2$. Then in the notation of the proof of
Theorem~\ref{lowerDurrett}, $\P (x \in A(t)) \leq e^ {nt} p(x, t)
\leq e^ {nt}
((1-e^ {-2t})/2)^ {n/2} ((1+e^ {-2t})/2)^ {n/2} = e^ {nt}
((1 - e^ {-4t})/4)^ {n/2} =
[(e^ {2t} - e^ {-2t}) / 4]^ n$. Now since $e^ {2t} < \sqrt{5}
+2$, it follows that
$e^ {2t} -e^ {-2t} < \sqrt{5}+ 2 - (\sqrt{5} - 2) = 4$, so $\P (x
\in A(t))
\rightarrow 0$ uniformly in such $x$ as $n \rightarrow \infty$.  Since at least
half the vertices of $\B_n$ satisfy $d(\0, x) \geq n/2$, it follows that
$\P(|A (t)^ c| \geq 2^ {n-2}) \rightarrow 1$.

Now condition on $|A(t)^ c| \geq 2^ {n-2}$.
The process $\{ A(s) : s \geq t \}$
is stochastically dominated by the process $\{ Z(u) : u \geq t \}$ for
which $Z(t) = A(t)$ and transitions from $S$ to $S \cup \{y\}$
occur at rate $n$ for all $S$ and $y \notin S$.  Now for each $x
\in A(t)^ c$, $\P (x \in Z(t+s)) = 1 - e^ {-ns}$, and
furthermore these events are independent as $x$ varies.  Thus
$\P (Z(t+s) = \B_n \| |A(t)^ c| \geq 2^ {n-2}) \leq
 (1-e^ {-ns})^ {2^ {n-2}}$.
Plugging in $s = \ln 2 - \ee / 2 < \ln(2-\ee /2)$ (for $\ee < 2$)
and using $1 - \ee < e^ {-\ee}$ gives
$$
\P (Z(t+s) = \B_n \| |A(t)^ c| \geq 2^ {n-2}) \leq
[e^ {-(2 - \ee/2)^ {-n}}]^ {2^ {n-2}} = e^ {-({2 \over 2 -
\ee / 2})^ n /
 4},
$$
which tends to $0$ as $n \rightarrow \infty$.  The theorem now follows
readily.  $\Cox$

\renewcommand{\baselinestretch}{1.0}\large

\begin{flushleft}
{\sc James Allen Fill \\
Department of Mathematical Sciences \\
The Johns Hopkins University \\
Baltimore, MD 21218-2689}\\
\end{flushleft}

%%%% I made this flushleft, since flushright looked dorky.
\begin{flushleft}
{\sc Robin Pemantle \\
Department of Mathematics, Van Vleck Hall \\
University of Wisconsin-Madison \\
480 Lincoln Drive \\
Madison, WI 53706}\\
\end{flushleft}

\end{document}